\documentclass[microtype]{gtpart}
\gtart % Mock-up which includes G&T header and footer.

\usepackage[all]{xy}

\title[Topological rigidity and involutions on tori]{Topological rigidity and\\ $H_1$-negative involutions on tori}

\author[F Connolly]{Frank Connolly}
\givenname{Frank}
\surname{Connolly}
\address{Department of Mathematics, University of Notre Dame, Notre Dame IN 46556 USA}
\email{connolly.1@nd.edu}
\urladdr{} % We'll leave this blank.

\author[J\,F Davis]{James F Davis}
\givenname{James F}
\surname{Davis}
\address{Department of Mathematics, Indiana University, Bloomington IN 47405 USA}
\email{jfdavis@indiana.edu}
\urladdr{} % We'll leave this blank.

\author[Q Khan]{Qayum Khan}
\givenname{Qayum}
\surname{Khan}
\address{Department of Mathematics, University of Notre Dame, Notre Dame IN 46556 USA}
\email{qkhan@indiana.edu}
\urladdr{} % We'll leave this blank.


\arxivreference{1102.2660} % [math.GT]
\arxivpassword{am5xz}

\volumenumber{}
\issuenumber{}
\publicationyear{}
\papernumber{}
\startpage{}
\endpage{}
\doi{}
\MR{}
\Zbl{}
\received{}
\revised{}
\accepted{}
\published{}
\publishedonline{}
\proposed{}
\seconded{}
\corresponding{}
\editor{}
\version{}

\newtheorem{thm}{Theorem}[section]
\newtheorem{conj}{Conjecture}[section]
\makeatletter \let\c@conj=\c@thm \makeatother
\newtheorem{cor}{Corollary}[section]
\makeatletter \let\c@cor=\c@thm \makeatother
\newtheorem{lem}{Lemma}[section]
\makeatletter \let\c@lem=\c@thm \makeatother
\newtheorem{prop}{Proposition}[section]
\makeatletter \let\c@prop=\c@thm \makeatother

\theoremstyle{definition}
\newtheorem{defn}{Definition}[section]
\makeatletter \let\c@defn=\c@thm \makeatother
\makeautorefname{defn}{Definition}
\newtheorem{rem}{Remark}[section]
\makeatletter \let\c@rem=\c@thm \makeatother
\newtheorem*{note}{Notation}

\DeclareMathAlphabet{\matheurm}{U}{eur}{m}{n}
\newcommand{\ab}{\matheurm{Ab}}
\newcommand{\all}{\matheurm{all}}
\newcommand{\fac}{\matheurm{fac}}
\newcommand{\fin}{\matheurm{fin}}
\newcommand{\Spectra}{\matheurm{Spectra}}
\newcommand{\topo}{\matheurm{Top}}
\newcommand{\vc}{\matheurm{vc}}

\newcommand{\RP}{\mathbb{RP}}
\newcommand{\F}{\mathbb{F}}

\newcommand{\bA}{\mathbf{A}}
\newcommand{\bE}{\mathbf{E}}
\newcommand{\bF}{\mathbf{F}}
\newcommand{\bG}{\mathbf{G}}
\newcommand{\bH}{\mathbf{H}}

\newcommand{\bL}{\mathbf{L}}

\newcommand{\bS}{\mathbf{S}}
\newcommand{\bT}{\mathbf{T}}

\newcommand{\mcA}{\mathcal{A}}
\newcommand{\mcB}{\mathcal{B}}
\newcommand{\mcC}{\mathcal{C}}
\newcommand{\cw}{\mathcal{CW}}
\newcommand{\mcD}{\mathcal{D}}
\newcommand{\mcF}{\mathcal{F}}
\newcommand{\mcG}{\mathcal{G}}

\newcommand{\mcP}{\mathcal{P}}

\newcommand{\mcS}{\mathcal{S}}
\newcommand{\mcT}{\mathcal{T}}
\newcommand{\mcU}{\mathcal{U}}

\newcommand{\tA}{\widetilde{A}}
\newcommand{\tf}{\widetilde{f}}
\newcommand{\tJ}{\widetilde{J}}
\newcommand{\tK}{\widetilde{K}}
\newcommand{\tM}{\widetilde{M}}
\newcommand{\tN}{{\widetilde{N}}}
\newcommand{\tX}{\widetilde{X}}
\newcommand{\tx}{\tilde{x}}
\newcommand{\G}{{\Gamma_n}}

\newcommand{\p}{{\partial}}

\newcommand{\bN}{\overline{N}}

\makeop{Aut}
\makeop{Cl}
\makeop{clc}
\makeop{cfd}
\makeop{Cok}
\makeop{Cyl}
\makeop{hAut}
\makeop{Ho}
\makeop{Holink}
\makeop{id}
\DeclareMathOperator{\Img}{Im}
\makeop{incl}
\makeop{Inn}
\makeop{Int}
\makeop{Ker}
\makeop{map}
\makeop{mic}
\DeclareMathOperator{\midd}{mid}
\makeop{Or}
\makeop{pt}
\makeop{res}
\makeop{UNil}
\makeop{Wh}

\newcommand{\TOP}{{\mathrm{TOP}}}

\newcommand{\bdry}{\partial}
\newcommand{\eps}{\varepsilon}
\newcommand{\longra}{\longrightarrow}
\newcommand{\ol}[1]{\overline{#1}}
\newcommand{\rel}{\;\mathrm{rel}\;}
\newcommand{\ul}[1]{\underline{#1}}
\newcommand{\connective}{\langle 1 \rangle}

\begin{document}

\begin{abstract}
We show, for $n\equiv 0,1 \pmod{4}$ or $n=2,3$, there is precisely one equivariant homeomorphism class of $C_2$-manifolds $(N^n, C_2)$ for which $N^n$ is homotopy equivalent to the $n$-torus and $C_2=\{1,\sigma\}$ acts so that $\sigma_*(x)=-x$ for all $x\in H_1(N)$. If $n\equiv 2,3 \pmod{4}$ and $n>3$, we show there are infinitely many such $C_2$-manifolds. Each is smoothable with exactly $2^n$ fixed points.  

The key technical point is that we compute, for all $n \geq 4$, the equivariant structure set $\mcS_\TOP(\R^n, \Gamma_n)$ for the corresponding crystallographic group $\Gamma_n$ in terms of the Cappell $\UNil$-groups arising from its infinite dihedral subgroups.
\end{abstract}

\maketitle

%------------------------------------------------------------------------------
%------------------------------------------------------------------------------
\section{Introduction}

%------------------------------------------------------------------------------
\subsection{Statement of results}\label{state}

Our goal here is to analyze   topological rigidity for a   sequence of crystallographic groups containing $2$-torsion. For each $n$, we define the group: $\Gamma_n = \Z^n \rtimes_{-1} C_2$, where $C_2$ acts on $\Z^n$ by negation: $v \mapsto -v$. 

We classify the proper actions of $\G$ on contractible $n$-manifolds.
   
The most powerful inspiration for our work is the remarkable  rigidity theorem of Farrell and Jones  concerning a discrete cocompact group of isometries of a simply connected non-positively curved manifold $(M, \Gamma)$. They classify the cocompact proper actions of such a $\Gamma$  on a contractible manifold, if $\Gamma$ is torsion-free.  

The second major  inspiration for our paper is the work of Cappell on $\UNil$.  If $\Gamma$ as above has elements of order 2, then the nontrivial elements of  $\UNil$ groups coming from virtually cyclic subgroups of $\Gamma$ can  provide  examples of cocompact $\Gamma$-manifolds $(M', \Gamma)$ which are isovariantly homotopy equivalent to, but not homeomorphic to $(M, \Gamma)$. So how do we classify such actions?

The   \emph{Topological Rigidity Conjecture} stated  below  does this. We view it as a   version of an old  conjecture of Quinn,  sharpened through the precision  afforded by the work of  \cite{ DL1}. We then prove  this conjecture  for  $\Gamma_n$ using \cite{BL_CAT0, CD}.

We can cast our results in terms of  an action of a  group $C_2:=\{1,\sigma\}$.  We say an involution $\sigma\co N\to N$ is \emph{$H_1$-negative} if $\sigma_*(x)=-x$, for all  $ x \in H_1(N)$. We prove:

\begin{thm}\label{thm:torus}
Let $\sigma\co N\to N$ be an $H_1$-negative involution on a closed manifold homotopy equivalent to the $n$-torus $T^n$.  Consider the $C_2$-manifold $(N, C_2)$.
\begin{enumerate}
\item\label{torus_discrete}
The fixed set $N^{C_2}$ is discrete and consists of exactly $2^n$ points.
\item\label{torus_rigid}
If $n \equiv 0,1\pmod{4}$ or $n=2,3$, then $(N^n, C_2)$ is equivariantly homeomorphic to the standard example, $(T^n, C_2)$.   
\item\label{torus_nonrigid}
If $n \equiv 2,3 \pmod{4}$ and $n > 3$, there are infinitely many such $C_2$-manifolds, $(N^n, C_2)$. All are isovariantly homotopy equivalent to $(T^n, C_2)$, but no two  are  equivariantly homeomorphic.  Each is smoothable hence locally linear.
\end{enumerate}
\end{thm}

By the \emph{standard example} $(T^n,C_2)$ above, we mean the involution $\sigma\co T^n \to T^n$ given by $\sigma[x] = [-x]$, for all $[x]\in \R^n/\Z^n = T^n$. Recall that any $n$-manifold homotopy equivalent to the $n$-torus is homeomorphic to it (see \cite{Wall_book, FQ, Anderson}).

The construction of the exotic involutions mentioned in the theorem uses surgery theory, specifically by Wall's realization \cite[Theorems~5.8, 6.5]{Wall_book}.  Write $X := (T^n - (T^n)^{C_2})/C_2$, an open $n$-manifold.  Define $\overline{X}$ as the obvious manifold compactification of $X$ obtained by adding a copy of $\RP^{n-1}$ at each end of $X$.  Note for all $n > 2$ that $\pi_1(\overline{X}) = \G$ and that $\overline{X}$ is orientable if and only if $n$ is even.  Let $w_n \co \G \to \{\pm 1\}$ be the orientation character of $\overline{X}$.  
Then, for $n \geq 5$, an element $\theta \in L_{n+1}(\G, w_n)$ determines a compact smooth manifold $\theta \cdot {\overline{X}}$, homotopy equivalent to $\overline{X}$ relative to the boundary.  Passing to the  two-fold cover and gluing in $2^n$ copies of $D^n$ with the antipodal action, we get a smooth involution on the torus.  All the exotic involutions in the above theorem arise in this way.

Observe that $\G$ is isomorphic to a rank $n$ crystallographic group. This isometric action of $\G$ on $\R^n$ is given by $\Z^n$ acting by translation and $C_2$ acting by reflection through the origin.  We let $(\R^n, \G)$ denote this $\G$-manifold.

Let $\mcS(\G)$ be the set of equivariant homeomorphism classes of contractible $n$-dimensional manifolds equipped with a proper $\G$-action.
We compute $\mcS(\G)$.

To parametrize the set $\mcS(\G)$, we will need to use the unitary nilpotent groups of Cappell. For $D_\infty$, these have been computed recently by Banagl, Connolly, Davis, Ko\'zniewski, and Ranicki \cite{CK_nilgps, CD, CR, BR}, yielding:
\begin{equation}\label{thm:BCDR}
\UNil_m(\Z;\Z,\Z) ~\cong~ \begin{cases}
0 & \text{if } m \equiv 0 \pmod{4}\\
0 & \text{if } m \equiv 1 \pmod{4}\\
(\Z/2\Z)^\infty & \text{if } m \equiv 2 \pmod{4}\\
 (\Z/2\Z \oplus \Z/4\Z)^\infty & \text{if } m \equiv 3 \pmod{4}.
\end{cases}
\end{equation}

Let $\midd(\Gamma_n)$ be the set of maximal infinite dihedral subgroups of $\Gamma_n$. Let $(\midd)(\Gamma_n)$ be a subset of $\midd(\Gamma_n)$ chosen so that it contains exactly one maximal infinite dihedral subgroup from each conjugacy class.   Let $D$ be a maximal infinite dihedral subgroup of $\G$.  For any integer $n$, with $\eps = (-1)^n$, there is a composite map
\[\xymatrix{
\UNil_{n+1}(\Z; \Z^\eps,\Z^\eps) \ar[r] & L_{n+1}(\Z D,w_n) \ar[r] & L_{n+1}(\Z\G,w_n).
}\]
If $n$ is odd, then there is an isomorphism $\UNil_{n-1}(\Z; \Z,\Z) \xrightarrow{\cong} \UNil_{n+1}(\Z; \Z^-, \Z^-)$.
 
\begin{thm}\label{thm:calculation}
Suppose $n \geq 4$. Write $\eps := (-1)^n$. The Wall realization map induces a  bijection of pointed sets, mapping the zero element to the basepoint $[\R^n,\Gamma_n]$: 
\[
\bdry_\oplus ~:~ \xymatrix{\displaystyle\bigoplus_{D \in (\midd)(\Gamma_n)} \UNil_{n+1}(\Z;\Z^\eps,\Z^\eps) \ar[r]^-{\approx} & \mcS(\G).}
\]
Consequently, $\mcS (\G)$ consists of a single element if $n \equiv 0, 1 \pmod{4}$, and $\mcS(\G)$ is countably infinite if  $n \equiv 2, 3 \pmod{4}$.
\end{thm}

We do not need to assume any conditions beyond continuity in order to obtain a full homeomorphism classification and to show all actions are smoothable.  It turns out that Smith theory guarantees the fixed sets consist of isolated points (see \fullref{sec:torus}). Also, local linearity is a consequence of our calculation (see \fullref{rem:nofakeboundary}, which concludes that the forgetful map $\mcS_\TOP(\ol{X}, \bdry\ol{X}) \to \mcS_\TOP(\bdry\ol{X})$ is constant).

An action $ \Gamma\times X \to X$ of a discrete group $\Gamma$ on a locally compact Hausdorff space $X$ is \emph{proper} if $\{\gamma\in \Gamma \mid K\cap \gamma K\neq \varnothing\}$ is finite for each compact set $K \subset X$.

Note that, given $(M,\G)\in \mcS(\G)$, the quotient manifold $M^n/\Z^n$ is homotopy equivalent to, and hence homeomorphic to the $n$-torus by the Borel Conjecture \cite{HW}. Therefore the universal cover $M$ admits a homeomorphism to $\R^n$.

\subsection*{Outline of the argument}
In \fullref{sec:torus}, we show that any $H_1$-negative involution on an $n$-manifold homotopy equivalent to the $n$-torus has exactly $2^n$ fixed points.  This allows one to deduce a correspondence between $H_1$-negative involutions on $n$-manifolds homotopy equivalent to the $n$-torus and contractible $n$-manifolds equipped with a proper $\G$-action.  In \fullref{sec:torus}
we show that any compact $C_2$-manifold with finite fixed set has the $C_2$-homotopy type of a finite $C_2$-CW-complex.  This allows one to conclude that any $H_1$-negative involution on a manifold homotopy equivalent to the $n$-torus is equivariantly homotopy equivalent to $(T^n,C_2)$ and that any  contractible $n$-manifold equipped with a proper $\G$-action is equivariantly homotopy equivalent to $(\R^n,\G)$.

We prove \fullref{thm:calculation} in \fullref{sec:calculation} and then deduce \fullref{thm:torus} in \fullref{sec:classify_involutions}.

For $n\geq 4$, the six structure sets we use are introduced in \fullref{isostruc}.  These are:
\[\begin{array}{ccccc}
\mcS(\G)
&&
\mcS_{\TOP}(\R^n,\G)
&&
\mcS^{iso}_{\TOP}(\R^n,\G)
\\ \\
\mcS_\TOP(\ol{X},\bdry\ol{X})
&&
\mcS_{\TOP}(T^n, C_2)
&&
\mcS^{iso}_{\TOP}(T^n, C_2).
\end{array}
\]
For example,  the isovariant structure set $\mcS^{iso}_{\TOP}(\R^n,\G)$ is the set of equivalence classes of proper $\G$-manifolds $(M^n, \G)$, together with an isovariant homotopy equivalence $M\to \R^n$.   We show all six structure sets are isomorphic, and compute the fourth one to prove \fullref{thm:calculation}.  The isomorphisms between the first and second, between the second and fifth, and the third and sixth structure sets are formal and are shown in \fullref{isostruc}.
The isomorphism between the fifth and sixth structure set requires a detailed discussion of equivariance versus isovariance and is discussed in \fullref{eq_to_iso}.  The isomorphism between the fourth and fifth structure set requires the use of end theory, see \fullref{lem:compactification}.  Finally, the computation of the classical surgery-theoretic structure set $\mcS_\TOP(\ol{X},\bdry\ol{X})$ uses the Farrell--Jones Conjecture and is presented at the end of \fullref{sec:calculation}.  This computation also uses the main result of \fullref{sec:QR=DL}, which identifies the assembly map in surgery theory with a corresponding map in equivariant homology.  
 
The final bit of the paper, \fullref{sec:quick}, is independent of \cite{BL_CAT0} and gives examples of non-standard structures on $(\R^n, \Gamma_n)$, hence of exotic $H_1$-negative involutions on tori.  The intent is to show that Cappell's work, for straightforward reasons, gives obstructions to isovariant rigidity of a $\Gamma$-space when $\Gamma$ has elements of order two.  Shmuel Weinberger pointed out these counterexamples to simple isovariant rigidity some time ago.  Since the argument was never published, we include it here.

\subsection{Equivariant Rigidity}

This paper represents the start of a systematic attack on Quinn's ICM conjecture and the closely related questions of equivariant, isovariant, and topological rigidity for a discrete group $\Gamma$.  We take some time to formulate these questions precisely.  Recall that two closed aspherical manifolds with the same fundamental group $\Gamma$ are homotopy equivalent and that the Borel Conjecture for $\Gamma$  states that any such homotopy equivalence is homotopic to a homeomorphism.  

Let $\Gamma$ be a discrete group.
A continuous function $f\co A \to B$ between $\Gamma$-spaces is \emph{equivariant}, or a \emph{$\Gamma$-map}, if $f(\gamma x) = \gamma f(x)$ for all $x \in A$ and $\gamma \in \Gamma$; it is \emph{isovariant} if furthermore $\gamma f(x)=f(x)$ implies $\gamma x=x$.
A \emph{model for $E_{\fin} \Gamma$} is a $\Gamma$-space $M$ $\Gamma$-homotopy equivalent to a $\Gamma$-CW-complex such that, for all subgroups $H$ of $\Gamma$,
\[
M^H \text{ is }
\begin{cases}
\text{contractible} & \text{if } H \text{ is finite}\\
\text{empty} & \text{otherwise}.
\end{cases}
\]
Given any $\Gamma$-CW-complex $X$ with finite isotopy groups, there is an equivariant map $X \to E_{\fin}\Gamma$, unique up to equivariant homotopy.  It follows that any two models are $\Gamma$-homotopy equivalent.  Furthermore, a model $E_{\fin}\Gamma$ exists for any group $\Gamma$.

A \emph{cocompact manifold model for $E_{\fin} \Gamma$} is a model $M$ for $E_{\fin} \Gamma$ so that $M/\Gamma$ is compact and so that $M^F$ is a manifold for all finite subgroups $F$ of $\Gamma$.  A geometric example is given by a discrete cocompact group $\Gamma$ of isometries of a simply connected complete nonpositively curved manifold $M$.
\emph{Equivariant (respectively, isovariant) rigidity holds for $\Gamma$} if any $\Gamma$-homotopy equivalence (respectively, $\Gamma$-isovariant homotopy equivalence) $M \to M'$ between cocompact manifold models for $E_{\fin} \Gamma$ is $\Gamma$-homotopic (respectively, $\Gamma$-isovariantly homotopic) to a homeomorphism.

With this terminology, our results can be restated as showing that every proper $\G$-action on a contractible manifold is a cocompact manifold model for $E_{\fin} \Gamma$, that equivariant and isovariant rigidity for $\G$ holds when $n \equiv 0,1 \pmod 4$ or $n = 2,3$, and that equivariant and isovariant rigidity fail for all other $n$.  Previous results on equivariant and isovariant rigidity are found in: \cite{Rosas_reflection}, \cite{CK_rigid}, \cite{CK_lack}, \cite[Section 14.2]{Weinberger_book}, \cite{Prassidis_Spieler}, and \cite{MP_rigidity}.  In particular, \cite{CK_lack} gave the first examples of groups where isovariant rigidity fails; this proceeded via a version of Whitehead torsion.  
\fullref{prop:Kflat} below shows that the relevant Whitehead group vanishes for $\G$.
We give the first counterexamples to simple isovariant rigidity in print.

In this paper we restrict ourselves to the study of equivariant and isovariant rigidity of $\G$, rather than for more general discrete groups $\Gamma$, for two reasons.   First, by  \fullref{thm:torus}\eqref{torus_discrete}, the singular set is discrete, and that simplifies the local analysis immensely.   Second, by \fullref{prop:Kflat} below, the group $\G$ is $K$-flat, so we avoid the subtleties of equivariant Whitehead torsion for topological manifolds.   In forthcoming work \cite{CDK2}, announced in \cite{Khan_announce},  we will study  equivariant rigidity of $E_{\fin}\Gamma$-manifolds with discrete singular set, without assuming the $K$-flat condition.

\begin{rem}\label{amalg}
A key algebraic property of $\G$ is that it admits  split epimorphisms
\[
\eps\co \G \longrightarrow \Gamma_1 = \Z \rtimes C_2 = C_2 * C_2
\]
to the infinite dihedral group.  The last equality results from noting that $\Gamma_1 = \langle (0,\sigma),(1,\sigma) \rangle = C_2 * C_2$. The existence of $\eps$ follows from the fact that every epimorphism $f \co \Z^n \to \Z$ gives a split epimorphism $f \rtimes \id \co \Z^n \rtimes C_2  \to \Z \rtimes C_2$.  Thus $\G$ has an injective amalgated product decomposition:
\[
\G ~=~ \eps^{-1}(C_2 * 1) *_{\eps^{-1}(1)} \eps^{-1}(1 * C_2) ~\cong~ \Gamma_{n-1} *_{\Z^{n-1}} \Gamma_{n-1}.
\]
\end{rem}

Our analysis of $E_\fin\G$ will have no issues with Whitehead torsion because:

\begin{prop}\label{prop:Kflat}
The group $\G$ is \emph{$K$-flat}, that is, $\Wh(\Gamma_n \times \Z^k) = 0$ for all $k \geq 0$.
\end{prop}

\begin{proof}
We prove $K$-flatness of $\G$ by induction on $n$, as follows.
First, note for $\Gamma_0 = C_2$ that $\Wh(C_2\times \Z^k)=0$, by using Rim's cartesian square of rings \cite[Section~3]{Milnor_Ktheory} and the vanishing of lower $NK$-groups of $\Z$ and $\F_2$ \cite[Chapter~XII]{Bass_Ktheory}.

Next, by \fullref{amalg} and Waldhausen's sequence \cite{Waldhausen_whitehead}, we obtain:
\[\xymatrix{
\Wh(\Gamma_{m-1}\times \Z^k)^{\oplus 2} \ar[r] & \Wh(\Gamma_m\times \Z^k) \ar[r]^-{\bdry} & \tK_0(\Z[\Z^{m-1+k}])=0.
}\]
This sequence is exact, since the $\mathrm{Nil}$ term vanishes because the ring $\Z[\Z^{m-1+k}]$ is regular coherent.
Therefore, by induction, we are done proving $\Gamma_n$ is $K$-flat.
\end{proof}

%------------------------------------------------------------------------------
\subsection{The Topological Rigidity Conjecture}\label{subsec:TRC}

This section is motivated by the conjecture of  F Quinn \cite{Quinn_ICM} at the 1986 ICM.
It aims to say the same thing, but in a more precise way, by employing the language of Davis--L\"uck \cite{DL1}.

Our \emph{Topological Rigidity Conjecture} concerns a discrete cocompact group $\Gamma$ of isometries of a simply connected complete nonpositively curved manifold $X^n$ (that is, a \emph{Hadamard manifold}).
It says, roughly, that any simple isovariant homotopy equivalence $f\co M\to X$ should be isovariantly homotopic to a homeomorphism, except for the examples created by $\UNil$-groups of virtually cyclic subgroups of $\Gamma$.
But it does so by parametrizing the set of such $f$ in terms of a homology group.
The coefficient spectrum of this homology group, $\ul{\bL}/\ul{\bL}_\fin$, is an $\Or(\Gamma)$-spectrum in the sense of  \cite{DL1} (see \fullref{sec:calculation}).
The homology is applied to a $\Gamma$-space with virtually cyclic isotropy.
For these virtually cyclic subgroups of $\Gamma$, the nonzero homotopy groups of the spectrum are just the $\UNil$ of amalgamated products of finite groups.

To formulate the conjecture, one must restrict to isovariant homotopy equivalences since there is no reason to expect equivariant homotopy equivalences to be well-behaved (see \cite[Section~14.2]{Weinberger_book} for some ill-behaved examples). We also restrict to \emph{simple} isovariant homotopy equivalences, whose definition is indicated below, to separate out the roles of $K$-theory and $L$-theory.   Furthermore,  \fullref{thm:torus}\eqref{torus_rigid} shows that this conjecture cannot be extended to low dimensions.

\begin{conj}\label{toprc}
Let $X^n$ be a Hadamard manifold of dimension $n>3$. Let $\Gamma$ be a discrete cocompact group of isometries of $X$.  Assume the fixed set $X^H$ has codimension $>2$ in $X^K$ whenever $K \subsetneqq H$ are isotropy groups.
There is a bijection 
\[\xymatrix{
H^\Gamma_{n+1}(E_{\vc}\Gamma; \ul{\bL}/\ul{\bL}_{\fin}) \ar[r]^-{\approx} & \mcS^{iso}_{rel}(X^n,\Gamma).
}\]
\end{conj}

The elements of $\mcS^{iso}_{rel}(X,\Gamma)$ are equivalence classes of pairs $(M,f)$, where
$M$ is a cocompact, locally flat topological $\Gamma$-manifold, and $f\co M \to X$ is a simple $\Gamma$-isovariant homotopy equivalence that restricts to a homeomorphism on the singular set.
\emph{Locally flat} means if $M^H \subset M^K$ then $M^H$ is a locally flat submanifold of $M^K$.
Here, $f\co M \to X$ and $f'\co M' \to X$ are equivalent if there is a $\Gamma$-homeomorphism $h\co M \to M'$ such that $f' \circ h$ is $\Gamma$-isovariantly homotopic to $f$.
Shmuel Weinberger has been a long-time proponent of this ``rel sing'' structure set in a very similar conjecture (for example, see \cite[Section~3]{CWY}).

If $\Gamma$ has no element of order 2, the conjecture implies each such $f$ is isovariantly homotopic to a homeomorphism.  The left side is defined in \cite{DQR}, using \cite{DL1}.

Quinn's conjecture should have included a vanishing  hypothesis on equivariant Whitehead torsion, as first noticed in \cite{CK_rigid, CK_lack}.
We state this simpleness condition using his subsequent work.
A $\Gamma$-isovariant homotopy equivalence $f\co M \to X$ between locally flat, cocompact $\Gamma$-manifolds is \emph{simple} if Quinn's stratified Whitehead torsion $\tau(\bar f \co M/\Gamma \to X/\Gamma)$ vanishes.
The element $\tau(\bar f)$ is defined using \cite[Corollary~1.6, Theorem~1.10(1)]{Quinn_HSS}.
The geometric interpretation of $\tau(\bar f)=0$ generalizes the interpretation for finite CW complexes and is given by \cite[Theorem~1.10(3)]{Quinn_HSS}: if $f$ is simple, then there is a homotopically stratified space $Z$ and a stratified cell-like map $z\co Z \to M/\Gamma$ such that $\bar f \circ z\co Z \to X/\Gamma$ is stratified homotopic to a cell-like map.

The Hadamard hypothesis is not, strictly speaking, necessary.   We could have stated the conjecture when $X$ is a locally flat, cocompact manifold model for $E_{\fin} \Gamma$ where each $X^H$ has codimension $>2$ in $X^K$.  We stated it in the Hadamard case to minimize jargon and to remind the reader that verification of \fullref{toprc} will require the Farrell--Jones Conjecture in $L$-theory for $\Gamma$, and that, to date, the verification of the Farrell--Jones Conjecture requires some geometric input.  In the more general case, we could conjecture,
\[\xymatrix{
H_{n+1}^\Gamma(X \longra \bullet; \ul{\bL}) \ar[r]^-{\approx} & \mcS^{iso}_{rel}(X^n,\Gamma).
}\]
This alternative conjecture would have the advantage that its verification should be independent of the Farrell--Jones Conjecture, but the disadvantage that the left hand side is not particularly computable without the Farrell--Jones Conjecture.

The present paper proves \fullref{toprc} for $\Gamma = \G$, using the facts that the singular set in $\R^n$ is discrete and that all finite subgroups of $\G$ have vanishing lower $K$-theory.
In fact our \fullref{thm:calculation} is stronger than \fullref{toprc} for $\Gamma_n$, and computes   $\mcS(\Gamma_n)$  using Smith theory and other topological tools. 
The bijection in \fullref{toprc} for $\Gamma=\G$ is defined in \fullref{sec:calculation} using maps of Cappell and Wall.

\subsection*{Acknowledgment}
The authors are grateful to Shmuel Weinberger  for suggesting that the work of Connolly--Davis \cite{CD} be exploited to construct counterexamples to an equivariant version of the Borel Conjecture,  as is done here. Years ago he mentioned the idea to the authors and has quietly encouraged the  project.  The second-named author thanks Mike Mandell for tutorials on the homotopy category.   

The National Science Foundation supports each author (DMS-0601234, 1210991, 0904276).  In addition this research has been facilitated by the support of the Hausdorff Institute, the Max Planck Institute, and the Humboldt Stiftung in Bonn.

%------------------------------------------------------------------------------
%------------------------------------------------------------------------------
\section{Applications of Smith theory}\label{sec:torus}

Given a universal covering map $p \co E \to B$ and an effective action of a group $G$ on $B$, consider the group
\[
\mcD (p,G) ~:=~ \{ h\in \mathrm{Homeo}(E) \mid p \circ h = g \circ p \text{ for some } g \in G \}.
\]
There is an obvious exact sequence
\begin{equation}\label{covgrpseq}\xymatrix{
1 \ar[r] & \mcD(p) \ar[r] & \mcD(p,G) \ar[r] & G \ar[r] & 1
}\end{equation} 
where $\mcD(p)$ is the deck transformation group of $p$ (those which cover the identity on $B$).

For the quotient map $q \co \R^n \to T^n$ and standard action $C_2 \curvearrowright T^n$, note $\mcD(q,C_2) = \Gamma_n$.

\begin{thm}\label{inv_on_torus}
Let $C_2$ act on a manifold $N^n$ homotopy equivalent to $T^n$  so that $\sigma_*(\alpha) = -\alpha$ for all $\alpha \in H_1(N^n)$.
\begin{enumerate}
\item\label{fixed_set}
The fixed set $N^{C_2}$ consists of $2^n$ points and $\mcD(p,C_2) \cong \Gamma_n$, where $p$ is the universal covering map of $N^n$. Moreover, if $G$ is any non-trivial finite subgroup of $\mcD(p,C_2)$, then ${\tN}^G$ consists of one point.

\item\label{iso_homotopy_unique}
Fix an isomorphism $\mcD(p,C_2) \cong \Gamma_n$. There is a $\Gamma_n$-homotopy equivalence  of the universal covers, $\tJ\co \tN^n \to \R^n$.   Any two such $\Gamma_n$-homotopy equivalences  are $\Gamma_n$-homotopic. $\tJ$ is the universal covering of a $C_2$-homotopy equivalence, $J\co N^n\to T^n$.
\end{enumerate}
\end{thm}

To prove this theorem we use  lemmas concerning involutions on $\R^n$ and $T^n$.  We state and prove these lemmas in their ultimate generality: actions of $p$-groups on contractible and aspherical manifolds.

\begin{lem}\label{Smith_theory}
Let $G$ be a finite $p$-group.
\begin{enumerate}
\item\label{Smith_theory_topman}
The fixed set of a $G$-action on a manifold is locally path-connected.
\item\label{Smith_theory_euclidean}
The fixed set of a $G$-action on a contractible manifold  is mod $p$ acyclic and path-connected.  
\item\label{Smith_theory_compact}
If the fixed set of a $G$-action on a contractible manifold is compact, then the fixed set is a point.
\end{enumerate}
\end{lem}

The proof of this lemma involves   Smith theory.  Our primary reference is Borel's Seminar on Transformation Groups \cite{Borel}.    Borel et.~al.~use Alexander--Spanier cohomology $\bar H^*(X; R)$ with coefficients in a commutative ring $R$.  This is, in turn, isomorphic to \v{C}ech cohomology $\check H^*(X; R)$ for $X$ paracompact Hausdorff, hence also for $X$ metrizable \cite[Corollary 6.8]{Spanier}.
Of course, if $X$ is a CW complex then Alexander--Spanier and \v{C}ech and singular cohomology coincide; but fixed sets of actions are far from CW complexes.    

The proof of this lemma is inductive; any non-trivial $p$-group has a normal $C_p$ subgroup, hence $G/C_p$ acts on $M^{C_p}$.  However the fixed set $M^{C_p}$ is not necessarily a manifold, so we will have to work abstractly.

\begin{proof}[Proof of \fullref{Smith_theory}]
\eqref{Smith_theory_topman}
Consider the following four properties of a topological space:  locally compact (lc), complete metrizable (cm), cohomology locally connected ($\clc_p$), and cohomology finite dimensional ($\cfd_p$).  We show that a manifold satisfies these four properties, that the fixed set of a $C_p$-action on a space which satisfies the four properties also satisfies the four properties, and that a space which satisfies the four properties is locally path-connected.  

A  space is {\em complete metrizable} if it admits a complete metric.  A space $X$ is \emph{cohomology locally connected with respect to $\F_p$} (written  $\clc_p$)  if each neighborhood $U$  of each $x\in X$ contains a neighborhood $V$ of $x$  so that the restriction map from $U$ to $V$ is zero on reduced \v{C}ech cohomology with coefficients in $\F_p$.  (See \cite[I.1.3]{Borel}  or \cite[II.17.1]{Bredon} for the definition.)  A locally compact space $X$ is  {\em cohomology finite dimensional} (written $\cfd_p$) if there is an integer $n \geq 0$ so that for all open sets $U$ of $X$, the $(n+1)$-st \v{C}ech cohomology with compact supports vanishes: $\check{H}^{n+1}_c(U;\F_p) = 0$.  (See \cite[I.1.2]{Borel} or \cite[II.16]{Bredon} for the definition.)  

Let $M$ be topological $n$-manifold.  Clearly it is locally compact.  We claim that it is complete metrizable.  Indeed by Urysohn's Metrization Theorem it admits the structure of a metric space $(M,d)$.  Since $M$ is locally compact and the union of a countable number of compact sets, it admits a proper map $f\co M \to \R$.  Then a complete metric on $M$ is given by $D(x,y) := d(x,y) + |f(x) - f(y) |$.  A manifold is $\clc_p$ since it is locally contractible.  
A topological $n$-manifold is $\cfd_p$; a reference is \cite[Corollary 16.28]{Bredon}; it also follows from Poincar\'e duality.

Thus it suffices to show that the fixed set of a $C_p$-action on a locally compact, locally path-connected, cohomology finite dimensional, complete metrizable space satisfies each of these properties.  Note that the fixed set $X^{C_p}$ is closed in $X$.  It clearly follows that $X^{C_p}$ is lc and cm, and it can be shown that any subspace of a $\cfd_p$ metrizable space is $\cfd_p$ (see \cite[Theorem~16.8]{Bredon}).  Finally, the fixed set of a $C_p$-action on a lc, $\cfd_p$, $\clc_p$ space is $\clc_p$ (see \cite[Proposition V.1.4]{Borel}, also \cite[Exercise~II.44]{Bredon}).  This is the one point where Smith theory is used.

Thus we conclude by induction on $|G|$ that the fixed set $M^G$ of a finite $p$-group acting on a manifold is $\clc_p$ and cm.  By looking at \v{C}ech cohomology in degree zero, it is easy to see that a $\clc_p$-space is locally connected (see \cite[page~6]{Borel} and \cite[page~126]{Bredon}).  A theorem of Moore--Menger--Mazurkiewicz \cite[page~254]{K} shows that a locally connected complete metric space is locally path-connected.

\noindent\eqref{Smith_theory_euclidean}
By a \emph{mod $p$ acyclic space}, we mean a Hausdorff space $X$ with $\check H^*(X; \F_p) \cong \check H^*(\pt; \F_p)$.  
A standard result from global Smith theory is that, if a finite $p$-group $G$  acts on a cohomology finite dimensional  mod $p$ acyclic space $X$, then its fixed set is also mod $p$ acyclic; see for example \cite[Cor.~III.4.6]{Borel}.  
In particular $\check H^0(M^{G}; \F_p) \cong \F_p$, so $M^{G}$ is connected.
 By Part \eqref{Smith_theory_topman}, it is locally path-connected.  Therefore $M^{G}$ is path-connected.

\noindent\eqref{Smith_theory_compact}
Since $M$ is a topological manifold and is contractible, $M$ is an orientable mod $p$ cohomology manifold  (see \cite[Definition 3.3]{Borel}.)  If $G$ acts on an orientable mod $p$  cohomology manifold $M^n$, then by Smith theory  the fixed set $M^G$ is an orientable mod $p$ cohomology manifold (\cite[Theorem V.2.2]{Borel}) of dimension $d \leq n$.  Recall  that $M^G$ is connected and mod $p$ acyclic.  If $M^G$ is compact, then there is a fundamental cohomology class in dimension $d$, so $\check H^d(M^G; \F_p) \cong \F_p$ (\cite[Theorem I.4.3(1)]{Borel}).  But $M^G$ is mod $p$ acyclic so $d = 0$.  Also, for a connected compact mod $p$ homology manifold $X$ of dimension $d$, $\check H^d(A; \F_p) =0$ for any closed proper subset $A$ of $X$ by \cite[Theorem I.4.3(1)]{Borel}.  Therefore $M^G$ is a point.
\end{proof}

In the literature, we have not seen the application of Smith theory to path-connectivity.  The above lemma seems to be the first occurrence.  We use path-connectivity in the covering space arguments below.

\begin{lem}\label{Smith}
Let $G$ be a discrete group.  Suppose  that $G$ acts effectively on a space $M$ fixing a  point $x \in M$.  Let $p \co (\tM, \tilde x) \to (M,x)$ be the universal cover.  
\begin{enumerate}

\item\label{Smith_lift} There is a unique action $G \curvearrowright (\tM, \tilde x)$ so that $p \co \tM \to M$ is a $G$-map.  

\item\label{Smith_inv} The action $G \curvearrowright (M,x)$ gives a homomorphism $G \to \Aut(\pi_1(M,x))$ and the action $G \curvearrowright  (\widetilde M, \tilde x)$ gives a splitting of the short exact sequence
\[\xymatrix{
1  \ar[r] & \mcD(p) \ar[r] & \mcD(p,G) \ar[r] & G \ar@<-1ex>@{-->}[l] \ar[r] & 1 .
}\]
The usual isomorphism $k \co \pi_1(M,x) \to \mcD(p)$ is equivariant with respect to the $G$-action on the fundamental group and the $G$-action given on $\mcD(p)$ given by conjugation.
Thus $k$ sends $\pi_1(M,x)^G$ to
\[
\mcD(p)^G = \{h \in \mcD(p) \mid ghg^{-1} = h \in \mcD(p,G) \text{ for all }g \in G\}.
\]

\item  \label{Smith_cover}  $p^G \co \widetilde M^G \to p(\widetilde M^G)$ is a regular $\mcD(p)^G$-cover.

\item  \label{Smith_path}  If $M^G_x$ is the path component of $M^G$ containing $x$, then $M^G_x \subset p(\widetilde M^G)$.  

\item  \label{Smith_contractible}  If $M$ is an aspherical manifold and $G$ is a finite $p$-group, then $p(\tM^G)$ is a connected component of $M^G$.

\item  \label{Smith_compact}  If $M$ is a compact, aspherical manifold and $G$ is a finite $p$-group, then $p(\tM^G)$ is compact.

\end{enumerate}
\end{lem}

\begin{proof}
\eqref{Smith_lift}   Any map $g\co (M,x) \to (M,  x)$ lifts to a unique map $\tilde g \co (\tM, \tilde x) \to (\tM, \tilde x)$ so that $p \circ \tilde g = g \circ p$.  Uniqueness implies that $\widetilde{g \circ h} = \tilde g \circ \tilde h$.  Define the action $G \curvearrowright \tM$ by $gy := \tilde g(y)$.

\noindent \eqref{Smith_inv} Clearly the $G$-action on $(M,x)$ gives a $G$-action on the fundamental group and the $G$-action on $(\tM,\tx)$ gives a splitting of the short exact sequence.    The isomorphism $k$ is specified by $k[\alpha](\tx) = \tilde \alpha (1)$ where $\tilde \alpha \co I \to \tM$ is the unique lift of $\alpha \co I \to M$ so that $\tilde \alpha (0) = \tx$. One then checks that $k$ is equivariant.

\noindent \eqref{Smith_cover}  Abstracting a bit, if $\pi\co E \to B$ is a regular $\Gamma$-cover, if $E_0$ is a closed subset of  $E$, and if $\Gamma_0$ is a subgroup of $\Gamma$ so that $\Gamma_0 \curvearrowright E_0$, and if $\gamma E_0$ and $E_0$ are disjoint whenever $\gamma \in \Gamma- \Gamma_0$, then $\pi_0 \co E_0 \to \pi(E_0)$ is a regular $\Gamma_0$-cover.  To check this in the case at hand, let $y \in \tM^G$.  If $h \in \mcD(p)^G$, then for all $g \in G$, $ghy = hgy = hy$, so $hy \in \tM^G$.  Conversely, if $h \in \mcD(p) - \mcD(p)^G$, then for some $g \in G$, $hg \not = gh$, so $ghy \not = hgy = hy$ by the freeness of the $\mcD(p)$-action, so $h y \not \in \tM^G$.

\noindent \eqref{Smith_path}  Let $y \in M^G_x$ and $\alpha \co I \to M^G$ be a path from $x$ to $y$.  Let $\tilde \alpha \co I \to \tM$ be the unique lift of $\alpha$ with $\tilde \alpha (0) = \tx$.  Then for any $g \in G$, $g \tilde \alpha$ is a lift starting at $g \tx = \tx$, so $g \tilde \alpha = \tilde \alpha$.  Hence $\tilde \alpha (1) \in \tM^G$ and since $p(\tilde \alpha(1) ) = y$ we have arrived at our conclusion.  

\noindent \eqref{Smith_contractible}  By \fullref{Smith_theory}\eqref{Smith_theory_euclidean}, $p(\tM^G)$ is path-connected and so by part \eqref{Smith_path}, $M^G_x = p(\tM^G)$.  On the other hand, by \fullref{Smith_theory}\eqref{Smith_theory_topman}, $M^G$ is locally path-connected so path components are components.  

\noindent \eqref{Smith_compact}  The fixed set $M^G$ is closed in $M$.  Connected components are closed, so $p(\tM^G)$ is closed in $M^G$.  Hence if $M$ is compact so is $p(\tM^G)$.
\end{proof}

\begin{proof}[Proof of \fullref{inv_on_torus}(1)]
Let $\sigma$ be an involution on a manifold $N^n$ homotopy equivalent to the $n$-torus such that $\sigma_* \alpha = - \alpha$  for all $\alpha \in H_1(N)$.   
  
We first note that $\sigma$ has a fixed point, because its Lefschetz number is $2^n\neq 0$.
Indeed, $L(\sigma) =  \sum_{k=0}^n (-1)^k \text{Tr}(\sigma_* \co H_k(N) \to H_k(N) )=  \sum_{k=0}^n (-1)^k(-1)^k\binom{k}{n}=2^n$.  Thus by \fullref{Smith}\eqref{Smith_inv}, $\mcD(p,C_2) \cong \Gamma_n$.  Since $\mcD(p) \cong \Z^n$ is torsion-free, all non-trivial finite subgroups of $\Gamma_n$ have order two and map isomorphically to $C_2$ under the map $\mcD(p,C_2) \to C_2$.

Let $G$ by a non-trivial finite subgroup of $\mcD(p,C_2)$.  We will show $\tN^G$ is a point.  By \fullref{Smith_theory}\eqref{Smith_theory_euclidean}, $\tN^G$ is nonempty.  Since $\mcD(p,C_2)^G$ is the trivial subgroup, by \fullref{Smith}\eqref{Smith_cover}, $\tN^G$ is homeomorphic to $p(\tN^G)$, which by \fullref{Smith}\eqref{Smith_compact} is compact.  Hence by \fullref{Smith_theory}\eqref{Smith_theory_compact}, $\tN^G$ is a point.

Next we must show $|N^{C_2}| = 2^n$. Each involution $s\in  \mcD(p,C_2)=\Gamma_n$ determines a point, $p(\tN^s)$ in $N^{C_2}$. Moreover,  for  involutions  $s$ and $s'$, we have: $p(\tN^{s})= p(\tN^{s'})$, iff for some deck transformation  $t \in \mcD(p), \;t(\tN^{s})=\tN^{s'}$, iff $ts't^{-1}=s$ for some $t\in \mcD(p)$, iff $s$ and $s'$ are conjugate in $\Gamma_n$ by an element $t\in \Z^n$, iff $s$ and $s'$ are conjugate in $\Gamma_n$ (because $\G = \Z^n\cup s\Z^n$).  Therefore this rule $s\mapsto p(\tN^s)$ induces a bijection between the set  of conjugacy classes of involutions in $\Gamma_n = \Z^n \rtimes C_2$ and $N^{C_2}$.  Since each such conjugacy class is represented uniquely by an element of $\{(\varepsilon, \sigma) \mid \varepsilon \in \{0,1\}^n\}$, we obtain $|N^{C_2}|=2^n$.
\end{proof}

Much of the above proof could be recast in a homological light.  The exact sequence \eqref{covgrpseq} is split, since $H^2(C_2; \mcD(p)) = 0$, so $\mcD(p,C_2) \cong \Gamma_n$. Then one shows, using the argument above, that for each non-trivial finite subgroup $G$ of $\mcD(p,C_2)$, $\tN^G$ is a point.    It follows that fixed points of $N^{C_2}$ are in one-to-one correspondence with conjugacy classes of splittings of \eqref{covgrpseq} which is in turn in one-to-one correspondence with $H^1(C_2; \mcD(p)) \cong (C_2)^n$.

Jiang \cite{Jiang} used fixed point theory to show that if $\tf \co \tX \to \tX$ is a lift of a self-map $f \co X \to X$ of a finite complex, then $p(\mathrm{Fix}~\tf)$ is compact.  This can be used to give an alternate proof that $p(\tN^{C_2})$, and hence $\tN^{C_2}$, is compact.

Before proving \fullref{inv_on_torus}(2), we need the following useful fact.

\begin{lem}\label{GCW}
Let $(N^n, C_2)$ be  compact $C_2$-manifold with finite fixed set.  If $n\leq 5$, assume   $\pi_1(N_{free}/C_2)\cong \G$.
Then $(N, C_2)$ has the equivariant homotopy type of a finite $C_2$-CW complex. 
\end{lem}

\begin{proof}
This is implicit, certainly in Quinn \cite{Quinn_HSS} but we will make it explicit here. First, the pair $(N/C_2,N^{C_2})$ is forward tame (tame in the sense of Quinn \cite{Quinn_HSS}) by \cite[Propositions 2.6 and 3.6]{Quinn_HSS}.  Second, the pair is reverse tame (tame in the sense of Siebenmann \cite{Siebenmann_thesis}) by \cite[Proposition 2.14]{Quinn_HSS}.   From this it follows from Siebenmann \cite{Siebenmann_thesis}, when $n\geq 6$,  that $N-N^{C_2}$ is the interior of a compact, free $C_2$ manifold, $\bN$, with $2^n$ boundary components and $N$ is homeomorphic to $\bN/\sim$, where $\sim $ is the equivalence relation identifying each boundary component to a separate point. But $\bN$ and each of its boundary components are homotopy equivalent to a finite free $C_2$ complex. It follows that $\bN/\sim $ is homotopy equivalent to a finite $C_2$ complex with $2^n$ fixed points.

Alternatively, in any dimension, we can argue that $(N-N^{C_2})/C_2$ and $\Holink(N/C_2, N^{C_2})$ are finitely dominated (see \cite[Proposition 2.15]{Quinn_HSS}). Since the projective class groups vanish: $\tK_0(\Z[C_2])=0=\tK_0(\Z[\G])$, by Wall's finiteness theorem~\cite[Thm.~F]{Wall_finiteness}, each is homotopy equivalent to a finite CW complex. So there is a finite $C_2$-CW pair $(K, L)$,  homotopy equivalent to the pair $(\Cyl (e_1), \Holink(N, N^{C_2}) )$, where $e_1\co \Holink(N, N^{C_2}) \to N-N^{C_2}$ is the evaluation map at time $1$. 
This map passes to  a $C_2$-homotopy equivalence of pairs:
\[\xymatrix{
(K\cup_L \Cyl(p), N^{C_2}) \ar[r] & (\Cyl (e_1)\cup_{H} \Cyl(e_0), N^{C_2})
}\]
where $H= \Holink(N, N^{C_2})$, and $p\co L\to N^{C_2}$ is the composite ${C_2}$-map, $L\to H \xrightarrow{e_0} N^{C_2}$. But, as noted in Quinn~\cite{Quinn_HSS},  $\Cyl (e_1)\cup_{H} \Cyl(e_0)$ is $C_2$-homotopy equivalent to $N$ and $K\cup_L \Cyl(p)$ is a finite $C_2$-CW complex.
\end{proof}

\begin{proof}[Proof of \fullref{inv_on_torus}(2)]
For a discrete group $\Gamma$, recall that a $\Gamma$-space $E$ is a \emph{model for $E_\fin\Gamma$} if   $E^G$ is contractible for each finite subgroup $G\subset \Gamma$,  and $E$ has the $\Gamma$-homotopy type of a $\Gamma$-CW complex, and  the   $\Gamma$-space $(E,\Gamma)$ is proper. 
Consider  $\R^n$ with its $\G$-action. The action is obviously proper. The only non-trivial finite subgroups have the form $\{1, \sigma\}$ for some involution $\sigma$; for this subgroup the fixed set is a single point. Finally, $\R^n$ admits the structure of a $\G$-CW complex. For it is the universal covering of $(T^n,C_2)$, and $(T^n,C_2)$ is the $n$-fold cartesian product of $(S^1, C_2)$, (with the diagonal action), which is a $C_2$-CW complex with exactly two (fixed) vertices.

Now let $(N, C_2)$ be as in the hypothesis of \fullref{inv_on_torus}. After the choice of an isomorphism $\mcD(p,C_2)\cong \G$, $\tN$ is a proper $\G$-manifold, and has the $\G$-homotopy type of a $\G$-CW complex by \fullref{GCW}. For each finite subgroup $G$ of $\G$, $\tN^G$ is contractible by \fullref{inv_on_torus}(1). Therefore $\tN)$ is also $\G$-universal, and there is a unique $\G$-homotopy class of $\G$-maps   $g\co \R^n\to \tN$. By uniqueness, $g$ and $\tJ$ are mutually $\G$-homotopy inverse. Therefore $\tJ \co   \tN \to \R^n$ and its quotient, $J\co N\to T^n$ are $C_2$-equivariant homotopy equivalences.
\end{proof}

This paper focuses on  actions of $C_2=\{1, \sigma\}$ on a torus $N$ for which $\sigma_*(x)=-x$ for all $x\in H_1(N)$.  But  the following lemma (with \fullref{inv_on_torus}) shows that this is \emph{equivalent} to saying that the torus has at least one isolated fixed point:

\begin{lem}
Suppose $C_2=\{1, \sigma\}$ acts on a torus $N$ and there is at least one fixed point which is isolated in $N^{C_2}$. Then $\sigma_*(x)=-x$ for all $x\in H_1(N)$.
\end{lem}

\begin{proof}
Lift the involution to an involution $\tilde{\sigma}\co \tN \to \tN$ with an isolated fixed point $\tx \in \tN$. The group $G=\{1,\tilde{\sigma}\}$ fixes only $\tx$ since the fixed set is mod 2 acyclic. But the centralizer, $Z_{\mcD(p)}(G)$ acts freely on this fixed set $\{\tx\}$, since the action is proper. So if $t\in \mcD(p)$, and $t \tilde{\sigma} = \tilde{\sigma} t$, then $t=1$.
Therefore, if $x\in H_1(N)\cong \mcD(p)$, and $ \sigma_*(x)=x$, then $x=0$. This implies that $\sigma_*(x)=-x$ for all $x\in H_1(N)$.
\end{proof}

%------------------------------------------------------------------------------
%------------------------------------------------------------------------------
\section{Equivariant and isovariant structures}\label{isostruc}

\begin{defn}\label{defn:isovariantstructureset}
The \emph {equivariant structure set}  $\mcS_\TOP(\R^n, \G)$,  is the set of equivalence classes of pairs $((M,\G), f)$, where $(M,\G)$ is a manifold with a cocompact proper $\Gamma_n$-action, and $f\co  M \to \R^n$ is a $\Gamma_n$-equivariant homotopy equivalence.  Often we write such pairs as $(M, f)$. Two such pairs $(M,f)$ and $(M', f')$ are equivalent if there is an equivariant homeomorphism $h\co  M \to M'$ and an equivariant homotopy $H$ from $f$ to $f'\circ h$. 
One defines the \emph{isovariant structure sets} $\mcS^{iso}_\TOP(\R^n,\G)$ and $\mcS^{iso}_\TOP(T^n,C_2)$ similarly, except one requires that $f$ is an isovariant homotopy equivalence and that $H$ is an isovariant homotopy. 
\end{defn}

\begin{rem}
We make no requirements on the fixed sets of subgroups in $M$ and $N$ above, because we have seen (in \fullref{sec:torus})  that  these fixed sets are discrete. A result of Quinn \cite[Prop.~2.6, Prop.~3.6]{Quinn_HSS} then guarantees that there are no local pathologies in such manifolds.
\end{rem}

The universal covering of a $C_2$-homotopy equivalence, $f\co  N^n  \to T^n$ is a $\Gamma_n$-homotopy equivalence $\tf\co  \tN^n \to \R^n$. Moreover $\tf$ is isovariant if $f$ is isovariant. This gives obvious bijections:
\begin{eqnarray}\label{iso:u}
u &:& \xymatrix{\mcS_\TOP(T^n, C_2) \ar[r]^-{\approx} & \mcS_\TOP(\R^n, \Gamma_n)}\\
u^{iso} &:& \xymatrix{\mcS^{iso}_\TOP(T^n, C_2) \ar[r]^-{\approx} & \mcS^{iso}_\TOP(\R^n, \Gamma_n).}
\end{eqnarray}

Consider the closed $C_2$-manifold $T:= T^n$.   Write $X := T_{free}/C_2$, an open $n$-manifold.  Define $\ol{X}$ as the obvious manifold compactification  of $X$ obtained by adding a copy of $\RP^{n-1}$ at each end of $X$. Label the boundary components:
\[
\bdry \ol{X} = \bigsqcup_{i=1}^{2^n} \bdry_i \ol{X}.
\]

\begin{defn}\label{defn:structureset_high} For $n\geq 5$, the \emph{structure set} $\mcS_\TOP(\ol{X},\bdry\ol{X})$ is the set of equivalence classes of triples $(\ol{Y},\ol{h},\bdry\ol{h})$, where $\ol{Y}$ is a compact $n$-manifold and $(\ol{h},\bdry\ol{h})\co  (\ol{Y},\bdry\ol{Y}) \to (\ol{X},\bdry\ol{X})$ is a homotopy equivalence of pairs.  Such triples $(\ol{Y}{}^0,\ol{h}{}^0,\bdry\ol{h}{}^0)$ and $(\ol{Y}{}^1,\ol{h}{}^1,\bdry\ol{h}{}^1)$ are equivalent if there is a homeomorphism $\varphi\co  \ol{Y}{}^0 \to \ol{Y}{}^1$ such that $(\ol{h}{}^1,\bdry\ol{h}{}^1) \circ \varphi$ is homotopic to $(\ol{h}{}^0,\bdry\ol{h}{}^0)$.
Compare \cite[\S10]{Wall_book}; we used $\Wh(\Gamma_n)=0$.
\end{defn}

For the four-dimensional case, some modifications are required.

\begin{defn}\label{defn:structureset_four}
For $n=4$, the \emph{structure set} $\mcS_\TOP(\ol{X},\bdry\ol{X})$ is the set of equivalence classes of triples $(\ol{Y},\ol{h},\bdry\ol{h})$, where $\ol{Y}$ is a compact topological 4-manifold and $(\ol{h},\bdry\ol{h})\co  (\ol{Y},\bdry\ol{Y}) \to (\ol{X},\bdry\ol{Y})$ is a $\Z[\Gamma_4]$-homology equivalence of pairs. Such triples $(\ol{Y}{}^0,\ol{h}{}^0,\bdry\ol{h}{}^0)$ and $(\ol{Y}{}^1,\ol{h}{}^1,\bdry\ol{h}{}^1)$ are equivalent if there is a $\Z[\Gamma_4]$-homology $h$-bordism $(\ol{W};\ol{Y}{}^0,\ol{Y}{}^1) \to \ol{X} \times (I;0,1)$ between them. Compare with \cite[\S11.3]{FQ}.
\end{defn}

For our application, we give a more explicit form of this general notion.

\begin{rem}\label{rem_ss4}
Since $\Wh(\Gamma_4)=0$ and $\Gamma_4$ is ``good'' in the sense of \cite{FQ}, we can simplify \fullref{defn:structureset_four}. First, for a representative $(\ol{Y},\ol{h},\bdry\ol{h})$, by the manifold-theoretic plus-construction $rel$ boundary \cite[\S11.1]{FQ}, we can assume that $\ol{h}\co  \bdry\ol{Y} \to \bdry\ol{X}$ is a homotopy equivalence. Second, two triples $(\ol{Y}{}^0,\ol{h}{}^0,\bdry\ol{h}{}^0)$ and $(\ol{Y}{}^1,\ol{h}{}^1,\bdry\ol{h}{}^1)$ are equivalent, by plus-construction on $\ol{W}$, if and only if there are:
\begin{itemize}
\item
a non-orientable closed 3-manifold $P = \bigsqcup_{i=1}^{16} P_i$

\item
$\Z[\Gamma_4]$-homology $h$-cobordisms $(\ol{Z}{}^j;\bdry\ol{Y}{}^j,P)$ for both $j=0,1$

\item
$\Z[\Gamma_4]$-homology equivalences $g^j\co  \ol{Z}{}^j \to \bdry\ol{X}$ extending $\bdry\ol{h}{}^j\co  \bdry\ol{Y}{}^j \to \bdry\ol{X}$

\item
a homeomorphism $\varphi\co  \ol{Y}{}^0 \cup_{\bdry\ol{Y}{}^0} \ol{Z}{}^0 \to \ol{Y}{}^1 \cup_{\bdry\ol{Y}{}^1} \ol{Z}{}^1$ relative to $P$
\end{itemize}
such that $(\ol{h}{}^1 \cup g^1) \circ \varphi$ is homotopic to $(\ol{h}{}^0 \cup g^0)$ relative to $P$. See \cite[\S11.1]{FQ}.
\end{rem}

\begin{lem}\label{lem:compactification}
Suppose $n \geq 4$. There is a bijection:
\[
\Phi ~:~ \xymatrix{\mcS_\TOP(\ol{X},\bdry\ol{X}) \ar[r]^-{\approx} & \mcS^{iso}_\TOP(T^n, C_2).}
\]
\end{lem}

\begin{proof}
First, suppose $n \geq 5$. Let $(\ol{h},\bdry\ol{h})\co  (\ol{Y}, \bdry\ol{Y}) \to (\ol{X}, \bdry\ol{X})$ be a homotopy equivalence of pairs, where $\ol{Y}$ is a compact $n$-dimensional topological manifold. So $[\ol{h},\bdry\ol{h}]$ is an element of $\mcS_\TOP(\ol{X}, \bdry\ol{X})$. Denote $\hat{X}$ and $\hat{Y}$ as the corresponding double cover of $\overline{X}$ and $\overline{Y}$.
Passage to double covers induces a $C_2$-equivariant homotopy equivalence $(\hat{h},\bdry \hat{h})\co  (\hat{Y}, \bdry\hat{Y}) \to (\hat{X}, \bdry\hat{X})$.
Each component $\bdry_i\hat{Y}$ of $\bdry\hat{Y}$ is homotopy equivalent to, and therefore homeomorphic to the sphere $S^{n-1}$ (see \cite{Smale_Poincare, FQ}). Hence the cone $c(\bdry_i\hat{Y})$  is homeomorphic to the disc $D^n$. So $N := \hat{Y} \cup \bigsqcup_i c(\bdry_i \hat{Y})$ is a topological $C_2$-manifold.
Thus, by coning off each map $\bdry_i \hat{h}$,  we obtain a function
\[
\Phi\co  \mcS_\TOP(\ol{X},\bdry\ol{X}) ~\longrightarrow~ \mcS^{iso}_\TOP(T, C_2) ~;~ (\ol{h},\bdry\ol{h}) ~\longmapsto~ \hat{h} \cup \bigsqcup_{i} c(\bdry_i \hat{h}).
\]
Now we show that $\Phi$ is a bijection by exhibiting its inverse.
Let $f\co  N \to T^n$ be a $C_2$-isovariant homotopy equivalence. This induces a proper homotopy equivalence
\[
f_{free}/C_2\co  N_{free}/C_2 ~\longrightarrow~ T_{free}/C_2.
\]
Then, since all the ends of $T_{free}/C_2$ are tame, the ends of $N_{free}/C_2$ are also tame.
Note $\Wh(C_2)=\tK_0(\Z[C_2])=0$.
Then, by a theorem of Siebenmann \cite{Siebenmann_thesis} (or of Freedman if $n=5$, see \cite{FQ}), we can fit a unique boundary, $\bdry\ol{N}$ onto $(N-N^{C_2})/C_2$, thereby creating a compact manifold $\ol{N}$, unique up to homeomorphism.
So we can extend  $f_{free}/C_2$ to $\bdry\ol{N}$.
(Here a small proper equivariant homotopy of $f_{free}/C_2$ may be needed before the extension.)
This construction, $[N, f]\mapsto [\ol{N}, \ol{f}]$ is clearly inverse to $\Phi$.
We conclude that $\Phi$ is both surjective and injective.   

Finally, it remains to consider $\Phi$ for $n=4$. Let $(\ol{h},\bdry\ol{h})\co  (\ol{Y}, \bdry\ol{Y}) \to (\ol{X}, \bdry\ol{X})$ be a map of pairs such that $\ol{h}\co  \ol{Y} \to \ol{X}$ is a homotopy equivalence of 4-manifolds and each $\bdry_i\ol{h}\co  \bdry_i\ol{Y} \to \bdry_i\ol{X}$ is a $\Z[C_2]$-homology equivalence of 3-manifolds.   Recall the notation $\hat{(~)}$ for the double cover of $\overline{(~)}$, used above for $n \geq 5$. For each $i$, by \cite[Proposition~11.1C]{FQ}, there is a compact contractible 4-manifold $c^*(\bdry_i\hat{Y})$ with a $C_2$-action such that its $C_2$-equivariant boundary is the homology 3-sphere $\bdry_i\hat{Y}$ and it has a single fixed point. It is unique up to $C_2$-homeomorphism. Using that isolated fixed point, one can construct a $C_2$-isovariant homotopy equivalence
\[
c^*(\bdry_i\hat{h}) ~:~ c^*(\bdry_i\hat{Y}) ~\longrightarrow~ c(\bdry_i\hat{X}).
\]
Suppose $(\ol{Y}{}^0,\ol{h}{}^0,\bdry\ol{h}{}^0)$ is equivalent to $(\ol{Y}{}^1,\ol{h}{}^1,\bdry\ol{h}{}^1)$, in the sense of \fullref{defn:structureset_four}.
In the setting of \fullref{rem_ss4}, there is a $C_2$-homeomorphism $\hat{Y}^0 \cup \hat{Z}^0  \to \hat{Y}^1 \cup \hat{Z}^1$. For each $j=0,1$, by uniqueness in \cite[Prop.~11.1C]{FQ}, there are $C_2$-homeomorphisms
\[
c^*(\bdry_i\hat{Y}^j) ~\longrightarrow~ \hat{Z}_i^j \cup c^*(\hat{P}_i).
\]
Using the identity map on each $c^*(\hat{P}_i)$, these produce a $C_2$-homeomorphism
\[
N^{0*} :=  \hat{Y}^0 \cup \bigsqcup_{i} c^*(\bdry_i\hat{Y}^0) ~\longrightarrow~ N^{1*} := \hat{Y}^1 \cup \bigsqcup_{i} c^*(\bdry_i\hat{Y}^1).
\]
Thus we may define $\Phi$ on equivalence classes similarly to the high-dimensional case, except we use the ``homotopy cones'' $c^*$ instead of the ``honest cones'' $c$. The argument for showing that $\Phi$ is a bijection as in the high-dimensional case, except for surjectivity we invoke the Weak End Theorem \cite[Thm.~11.9B]{FQ} and for injectivity we invoke the Classification of Weak Collars \cite[Thm.~11.9C(3)]{FQ}.
\end{proof}

\begin{lem}\label{lem:nobad}
Suppose $n \geq 4$. Consider the above-defined structure sets.
\begin{enumerate}
\item
The following forgetful maps are bijections:
\begin{eqnarray*}
\Psi &:& \xymatrix{\mcS^{iso}_\TOP(\R^n,\G) \ar[r]^-{\approx} & \mcS_\TOP(\R^n,\G)}\\
\psi &:& \xymatrix{\mcS^{iso}_\TOP(T^n,C_2) \ar[r]^-{\approx} & \mcS_\TOP(T^n,C_2).}
\end{eqnarray*}
\item
The following forgetful map is a bijection:
\begin{eqnarray*}
\chi &:& \xymatrix{\mcS_{\TOP}(\R^n,\G) \ar[r]^-{\approx} & \mcS(\G).}
\end{eqnarray*}
\end{enumerate} 
\end{lem}

\begin{proof}
For Part~(1), in view of the bijections $u$ and $u^{iso}$, it suffices to prove that $\psi$ is a bijection.
This is achieved by \fullref{equi-iso}, which we provide in \fullref{eq_to_iso}.

For Part~(2), by \fullref{inv_on_torus}(2), it is immediate that $\chi$ is injective.
We must prove $\chi$ is surjective.
Let $[M,\G]$ be in $\mcS(\G)$.
Since $M$ is a contractible $n$-manifold and the restricted action $\Z^n  \curvearrowright M$ is free, the quotient space $M/\Z^n$ is an $n$-manifold and an Eilenberg--MacLane space of form $ K(\Z^n, 1)\simeq T^n$.
The action $\G \curvearrowright M $ descends to $C_2  \curvearrowright M/\Z^n  $ with $(\alpha \mapsto -\alpha)$ on $H_1$.
So, by \fullref{inv_on_torus}(2), there is a $\G$-homotopy equivalence $\tJ\co  M \to \R^n$, yielding an element $[(M,\G), \tJ]\in \mcS_{\TOP}(\R^n, \G)$.
Thus $\chi$ is surjective.
\end{proof}

\begin{prop}\label{prop:alpha}
Suppose $n \geq 4$.
The following map $\alpha$ is a bijection:
\[
\alpha ~:=~ \chi \circ u \circ \psi \circ \Phi ~:~ \xymatrix{\mcS_\TOP(\ol{X},\bdry\ol{X}) \ar[r]^-{\approx} & \mcS(\Gamma_n)} ~;~\quad [\ol{X},\id] \longmapsto [\R^n,\Gamma_n].
\]
\end{prop}

\begin{proof}
This follows immediately from \eqref{iso:u}, \fullref{lem:compactification}, and \fullref{lem:nobad}.
\end{proof}

%------------------------------------------------------------------------------
%------------------------------------------------------------------------------
\section{Calculation of the isovariant structure set}\label{sec:calculation}

Our ultimate goal here is to prove \fullref{thm:calculation}. We also establish the Topological Rigidity Conjecture (of \fullref{subsec:TRC}) for the crystallographic groups $\Gamma_n$.

Throughout this section, we assume $n\geq 4$ and shall use the shorthand $\Gamma := \Gamma_n$. For each family $\mathcal{F}$ of subgroups of $\Gamma$, we write $E_{\mathcal{F}}\Gamma$ for the classifying space for $\Gamma$-CW complexes whose isotropy groups are in $\mathcal{F}$. 
We use the families $\fin$, $\vc$, and $\all$, consisting  of finite, virtually cyclic  and all subgroups respectively. For the remainder of this section, since the subgroups of $\Gamma=\Gamma_n$ have trivial reduced lower $K$-theory, for ease of reading, we shall simply write $\ul{\bL}$ for the $\Or(\Gamma)$-spectrum $\ul{\bL}^h$.

Recall the Wall realization map \cite[Thms.~5.8, 6.5]{Wall_book}, relative to the boundary:
\[
\bdry^{Wall} ~:~ L_{n+1}^h(\Gamma,w_n) \longrightarrow \mcS_\TOP(\ol{X}) \longrightarrow \mcS_\TOP(\ol{X},\bdry\ol{X}).
\]

Using Cappell's map \cite{Cappell_unitary}, define a composite homomorphism
\[
\beta ~:~ \bigoplus_{D \in (\midd)(\Gamma)} \UNil_{n+1}(\Z;\Z^\eps,\Z^\eps)
\longrightarrow \bigoplus_{D \in (\midd)(\Gamma)} L_{n+1}^h(D,w_n)
\longrightarrow L_{n+1}^h(\Gamma,w_n).
\]

Now we can define the desired basepoint-preserving function
\[
\bdry_\oplus ~:=~ \alpha \circ \bdry^{Wall} \circ \beta ~:~ \bigoplus_{D \in (\midd)(\Gamma)} \UNil_{n+1}(\Z;\Z^\eps,\Z^\eps) \longrightarrow \mcS(\Gamma).
\]

It remains to show that $\bdry_\oplus$ is a bijection of sets. This will span several lemmas.

%------------------------------------------------------------------------------
\subsection{Algebraic structure groups and equivariant homology}

For cofibrant pairs $(A,B)$ of topological spaces, A.~Ranicki \cite{Ranicki_TSO} defined the algebraic structure groups $\mcS^h_*$ as the homotopy groups of the homotopy cofiber of an assembly map $\alpha\connective$, so that there is a long exact sequence 
\[
\cdots \to H_*(A,B;\bL\connective) \xrightarrow{\alpha\connective} L_*^h(A,B) \to \mcS^h_*(A,B) \to H_{*-1}(A,B;\bL\connective) \xrightarrow{\alpha\connective} \cdots.
\]
Here $\bL\connective$ is the $1$-connective cover of the $4$-periodic surgery spectrum $\bL$ algebraically defined in \cite{Ranicki_TSO}. 
(The homotopy invariant functor $\mcS^h_*$ is a desuspended chain-complex analogue of the geometric structure groups $\mcS^h_*$ originally defined by F~Quinn.) When a map $i\co  B \to A$ is understood, we shall write $(A,B)$ for the cofibrant pair $(\Cyl(i),B)$. The relative $L$-groups $L^h_*(A,B) = L^h_*(i\co  B \to A)$ were defined algebraically by Ranicki (\cite{Ranicki_ESATS}), following C\,T\,C~Wall \cite{Wall_book}.

For computational purposes, we employ the non-connective version $\mcS^{per,h}_*$ of $\mcS^h_*$. It is the homotopy groups of a homotopy cofiber of an assembly map $\alpha$:
\[
\cdots \to H_*(A,B;\bL) \xrightarrow{\alpha} L_*^h(A,B) \to \mcS^{per,h}_*(A,B) \to H_{*-1}(A,B;\bL) \xrightarrow{\alpha} \cdots.
\]

\begin{rem}\label{rem:assemblymaps_same}
Suppose $B = \varnothing$ and that $A$ is the quotient of a free $\Gamma$-action on a space $\tA$ each of whose components is simply connected. Write $\Pi_0(\tA)$ for the $\Gamma$-set of components of $\tA$; there is a canonical $\Gamma$-map $\tA \to \Pi_0(\tA)$.
By \fullref{QR=DL}, the Quinn--Ranicki assembly map can be naturally identified with the Davis--L\"uck assembly map, at the spectrum level.  Then the cofibers of these assembly maps agree in a functorial manner. Specifically, \fullref{sec:QR=DL} constructs an isomorphism in $\Ho \mathrm{sc}(\Gamma,1)\cw\text{-}\Spectra$, whose value on $\tA$ after the application of homotopy groups gives an isomorphism:
\begin{equation}\label{eqn:relativeBredon}
\xymatrix{H^{\Gamma}_*(\Pi_0(\tA), \tA; \ul{\bL}) \ar[r]^-{\cong} & \mcS_*^{per,h}(A).}
\end{equation}
\end{rem}

Write $\R^n_{free} := \{ x \in \R^n \mid \Gamma_x = 1 \}$ for those points with trivial isotropy group.
Observe that $\R^n_{free}$ equivariantly deformation retracts to the universal cover of $\ol{X}$.
There is a canonical $\Gamma$-map from $\R^n$ to its singleton $\{\R^n\}$ with trivial $\Gamma$-action.

\begin{lem}\label{lem:split}
There is a commutative diagram with long exact rows:
\[\small\xymatrix{
H_*^\Gamma(\R^n,\R^n_{free}; \ul{\bL}) \ar[r] \ar[d] & H_*^\Gamma(\{\R^n\},\R^n_{free};\ul{\bL}) \ar[r] \ar[d] & H_*^\Gamma(\{\R^n\}, \R^n; \ul{\bL}) \ar[r] \ar[d] & H_{*-1}^\Gamma(\R^n,\R^n_{free}; \ul{\bL}) \ar[d]\\
\mcS^{per,h}_*(\bdry\ol{X}) \ar[r] & \mcS^{per,h}_*(\ol{X}) \ar[r] & \mcS^{per,h}_*(\ol{X},\bdry\ol{X}) \ar[r] & \mcS^{per,h}_{*-1}(\bdry\ol{X}).
}\]
Furthermore, the vertical maps are isomorphisms of abelian groups.
\end{lem}

\begin{proof}
The top row is the long exact sequence of the triple $(\{\R^n\}, \R^n, \R^n_{free})$ in $\Gamma$-equivariant $\bL$-homology \cite{DL1}.
The bottom row is the long exact sequence of the pair $(\ol{X},\bdry\ol{X})$ in algebraic structure groups \cite{Ranicki_TOP}.
The inner left vertical map is induced on homotopy groups from \eqref{eqn:relativeBredon} for $\tA = \R^n_{free}$. Write $D(\R^n_{sing})$ for the $\Gamma$-subset of points in $\R^n$ of distance $\leq 0.2$ from $\R^n_{sing}$. Define $S(\R^n_{sing}) := \bdry D(\R^n_{sing})$. The outer vertical maps are induced on homotopy groups from \eqref{eqn:relativeBredon} for $\tA = S(\R^n_{sing})$, precomposed with the inverse of the induced excision map
\[
H_*^\Gamma(D(\R^n_{sing}), S(\R^n_{sing}); \ul{\bL}) \longra H_*^\Gamma(\R^n, \R^n_{free}; \ul{\bL}).
\]
By functoriality of \eqref{eqn:relativeBredon}, the left square's diagram of spectra homotopy-commutes. In particular, the left square itself commutes. Then the inner right map is induced from a well-defined homotopy class of map of spectra.
So the middle square and right square are defined and commute. Therefore, since \eqref{eqn:relativeBredon} implies the outer maps and inner left map are isomorphisms, by the Five Lemma, we conclude that the inner right map is an isomorphism also.
\end{proof}

\begin{rem}\label{rem:nofakeboundary}
By the Isomorphism Conjecture \cite{BL_CAT0} and Bartels' splitting theorem \cite{Bartels_splitting} on the top row, we conclude that the connecting homomorphism $\mcS^{per,h}_*(\ol{X},\bdry\ol{X}) \to \mcS^{per,h}_{*-1}(\bdry\ol{X})$ is zero. Therefore, using Ranicki's natural bijection (see \fullref{rem:TSO} below), the forgetful map $\mcS_\TOP(\ol{X},\bdry\ol{X}) \to \mcS_\TOP(\bdry\ol{X})$ is constant.
\end{rem}

We now calculate $H^\Gamma_*(E_\vc\Gamma,E_\fin\Gamma; \ul{\bL})$ by using a specific model for the spaces involved. Models of $E_\vc G$ for crystallographic groups $G$ are due to Connolly--Fehrman--Hartglass \cite{CFH}. For any group $G$, L\"uck--Weiermann \cite{LW} built models of $E_\vc G$ from $E_\fin G$. However, the following lemma is shown directly for our $\Gamma$.

Let $C$ be an infinite cyclic subgroup of $\Gamma$. Let $\mcP_C$ denote the collection of all affine lines $\ell \subset \R^n$ which are stablized by $C$. Endow $\mcP_C$ with the affine structure and topology of a copy of $\R^{n-1}$. Since $\mcP_C$ is a partition of $\R^n$, there is a continuous quotient map $\pi_C\co  \R^n \to \mcP_C$. Since $C$ is normal in $\Gamma$, the $\Gamma$-action on $\R^n$ extends to a $\Gamma$-action on the mapping cylinder, $\Cyl(\pi_C)$.
  
Let $\mic(\Gamma)$ denote the collection of maximal infinite cyclic subgroups of $\Gamma$.

\begin{lem}\label{lem:model}
A  model $E$ for the classifying space $E_\vc\Gamma$ (classifying $\Gamma$-CW complexes with virtually cyclic isotropy) is the union along $\R^n$ of mapping cylinders:
\[
E := \bigcup_{C \in \mic(\Gamma)} \Cyl(\pi_C\co  \R^n \longrightarrow \mcP_C).
\]
\end{lem}

\begin{proof}
If $H$ is a finite nontrivial subgroup of $\Gamma$, then $E^H$ is a tree with one edge in $\Cyl(\pi)^C$ for each $C\in \mic(\Gamma)$. So $E^H$ is contractible. If $H$ is infinite cyclic or infinite dihedral, there is just one $C\in \mic(\Gamma)$ for which $\Cyl(\pi_C)^H$ is nonempty. For this $C$, observe that $\Cyl(\pi_C)^H$ is a single point when $H$ is dihedral and is all of $\Cyl(\pi_C)$ when $H$ is cyclic. Also $E = E^{\{1\}}$ is contractible, and $E^H$ is empty if $H$ is not virtually cyclic.

Finally, we  must prove that $E$ has the structure of a $\Gamma$-CW complex. We begin by assuming $K$ is a  $\Gamma$-CW structure on $\R^n$ which is \emph{convex}. By this we mean each closed cell is convex, and its boundary is a subcomplex.  It suffices to show how to extend $K$ to a $\Gamma$-CW structure over each mapping cylinder, $\Cyl(\pi_C)$ in $E$.

So fix $C$ and parametrize $\Cyl(\pi_C)$ as 
\[
\Cyl(\pi_C)= \R^n\times[-1,1]\cup_{\pi_C}\mcP_C , \qquad \text{where $(x,1)=\pi_C(x)$ for all $x\in \R^n$}.
\]

There are convex $\Gamma$-CW structures, $L$ on $\mcP_C$, and $\hat{L}$ on $\R^n$, so that  each $j$-cell $f$ of $L$ has the form $\pi_C(\hat{f})$ for some $(j+1)$-cell $\hat{f}$ of $\hat{L}$. This endows $\R^n\times [0, 1]\cup_{\pi_C}\mcP_C$    with the structure of a $\Gamma$-CW complex, $K_+$ so that $\R^n\times 0$ is the complex $\hat{L}$. Now $\hat{L}$ and $K$ have a common subdivision $K'$,  since each is convex. There is then a CW structure $K_-$ on $\R^n\times [-1,0]$ in which $K$, $K'$ and $\hat{L}$ are identified with $\R^n\times \{-1\}, \R^n\times \{-\frac{1}{2}\} $ and $ \R^n\times \{0\} $ respectively as subcomplexes. (Also, $e\times[-1,-\frac{1}{2}]$ and $f\times [-\frac{1}{2}, 0]$ are cells of $K_-$  if $e$ and $f$ are cells of $K$ and $\hat{L}$ respectively.) Then $K_+\cup K_-$ is the required $\Gamma$-CW structure on $\Cyl(\pi_C)$.
\end{proof}

Each infinite dihedral subgroup $D$ of $\Gamma$ contains a unique maximal infinite cyclic subgroup $C$. Moreover, $D$ has a unique invariant line, $\ell_D \subset \R^n$.
The image of $\ell_D$ in $\mcP_C$ is a single point, which we denote by the singleton $\{\ell_D\} = \pi_C(\ell_D)$.

\begin{lem}\label{lem:uncontrolledcalculation}
The inclusion-induced map is an isomorphism of abelian groups:
\[
\bigoplus_{D \in (\midd)(\Gamma)} H^D_*(\{\ell_D\}, \ell_D; \ul{\bL}) ~\longrightarrow~ H_*^\Gamma(E,\R^n; \ul{\bL}).
\]
\end{lem}

\begin{proof}
Lemma 4.1 of \cite{DQR} allows one to translate between maps induced by $\Gamma$-maps of classifying spaces for actions with isotropy in a family to maps induced by maps of $ \Or(\Gamma)$-spectra.  There is a homotopy cofiber sequence of  $\Or(\Gamma)$-spectra:
\[
\bL_{\fin} \longrightarrow \bL \longrightarrow \bL/\bL_{\fin}.
\]
By \cite[Lemma~4.1(ii)]{DQR}, the following absolute homology group vanishes:
\[
H_*^\Gamma(\R^n;  \ul{\bL} / \ul{\bL}_\fin ) ~=~ 0.
\]
Also, by \cite[Lemma~4.1(iii)]{DQR}, the following relative homology group vanishes:
\[
H_*^\Gamma(E,\R^n;\ul{\bL}_\fin) ~=~ 0.
\]
So we obtain a composite isomorphism, informally first observed by Quinn:
\begin{equation}\label{iso:Quinn}
\xymatrix{H_*^\Gamma(E,\R^n; \ul{\bL}) \ar[r]^-{\cong} & H_*^\Gamma(E,\R^n; \ul{\bL} / \ul{\bL}_\fin)
& \ar[l]_-{\cong} H_*^\Gamma(E; \ul{\bL}/\ul{\bL}_\fin).}
\end{equation}
Now, since $N_\Gamma(C)=\Gamma$, by \fullref{lem:model} and excision, we obtain:
\begin{equation}\label{evclolfin}
\xymatrix{\displaystyle\bigoplus_{C \in \mic(\Gamma)} H_*^\Gamma(\mcP_C; \ul{\bL}/\ul{\bL}_\fin)
\ar[r]^-{\cong} & H_*^\Gamma(E; \ul{\bL}/\ul{\bL}_\fin).}
\end{equation}
Fix $C\in \mic(\Gamma)$.  Observe that the action of the group $\Gamma/C$ on the parallel pencil $\mcP_C$ has a discrete singular set:
\[
{sing}\,\mcP_C ~:=~ \{ \ell_D \in \mcP_C \mid C \subset D \text{ for some } D\in \midd(\Gamma) \}.
\]
Let $U$ be a $\Gamma$-tubular neighborhood of ${sing}\,\mcP_C$ in $\mcP_C$. Write $V := \mcP_C - {sing}\,\mcP_C$.
Recall, by a theorem of J.~Shaneson \cite{Shaneson}, that the following assembly map is a homotopy equivalence:
\[
A_C ~:~ \xymatrix{S^1_+ \wedge \bL(1) \ar[r]^-{\simeq} & \bL(C).}
\]
That is, the spectrum $(\ul{\bL}/\ul{\bL}_\fin)(\Gamma/C)$ is contractible. So, since $V$ has isotropy $C$:
\[
H^\Gamma_*(U\cap V; \ul{\bL}/\ul{\bL}_\fin) ~=~ 0 ~=~ H^\Gamma_*(V; \ul{\bL}/\ul{\bL}_\fin).
\]
Since ${sing}\,\mcP_C$ is discrete, there is a $\Gamma$-homotopy equivalence:
\[\xymatrix{
\displaystyle\bigsqcup_{\substack{D \supset C\\ D \in (\midd)(\Gamma)}} \Gamma \times_D \{\ell_D\}
\ar[r]^-{\cong} & {sing}\,\mcP_C
\ar[r]^-{\simeq} & U.
}\]
So the homotopy and excision axioms of equivariant homology imply:
\[\xymatrix{
\displaystyle\bigoplus_{\substack{D \supset C\\ D \in (\midd)(\Gamma)}} H_*^\Gamma(\Gamma \times_D\{\ell_D\}; \ul{\bL}/\ul{\bL}_\fin)
\ar[r]^-{\approx} &
H_*^\Gamma(U; \ul{\bL}/\ul{\bL}_\fin)
\ar[r]^-{\cong} &
H_*^\Gamma(\mcP_C; \ul{\bL}/\ul{\bL}_\fin).
}\]
Thus \eqref{evclolfin} and the induction axiom of equivariant homology imply:
\[\xymatrix{
\displaystyle\bigoplus_{ D \in (\midd)(\Gamma)} H_*^D( \{\ell_D\} ; \ul{\bL}/\ul{\bL}_\fin)
\ar[r]^-{\cong} & H_*^\Gamma(E; \ul{\bL}/\ul{\bL}_\fin). 
}\]
Finally, since $\ell_D$ is a model for $E_\fin D$, by \cite[Lemma~4.1(ii)]{DQR} again, we obtain:
\[\xymatrix{
H^D_*( \{\ell_D\},\ell_D; \ul{\bL})
\ar[r]^-{\cong} &
H_*^D( \{\ell_D\}, \ell_D ; \ul{\bL}/\ul{\bL}_\fin)
& \ar[l]_-{\cong}
H_*^D( \{\ell_D\} ; \ul{\bL}/\ul{\bL}_\fin).
}\proved\]
\end{proof}

\begin{lem}\label{lem:UNil_relhom}
Recall $\eps = (-1)^n$. Let $D$ be an infinite dihedral subgroup of $\Gamma$. Then the following composite map is an isomorphism of abelian groups:
\[
\UNil_{n+1}(\Z;\Z^\eps,\Z^\eps) \longrightarrow L_{n+1}(D,w_n) = H_{n+1}^D(\{\ell_D\}; \ul{\bL}) \longrightarrow H_{n+1}^D(\{\ell_D\}, \ell_D; \ul{\bL}).
\]
\end{lem}

\begin{proof}
Denote the map under consideration by $\phi$. Consider the three maps
\begin{gather*}
\xymatrix{
L_{n+1}(C_2,\eps)\oplus L_{n+1}(C_2,\eps) \ar[r]^-{i} & H_{n+1}^D(\ell_D;\ul{\bL}) \ar[r]^-{j} & L_{n+1}(D,w_n)
}\\
\tag*{\text{and}} 
\xymatrix{
L_n(1) \ar[r]^-{k} & L_{n}(C_2,\eps)\oplus L_{n}(C_2,\eps).
}
\end{gather*}
Observe that $\phi$ has a factorization given by the commutative diagram
\[\xymatrix{
& \Cok(i) \ar[d]_-{j_*} \ar[r]^-{\bdry^{alg}} & \Ker(k)\\
\UNil_{n+1}(\Z;\Z^\eps,\Z^\eps) \ar[d]_-{\phi} \ar[r]^-{c} & \Cok(j \circ i) \ar[d]_-{\kappa} \ar[r]^-{\bdry^{top}} & \Ker(k)\\
H_{n+1}^D(\{\ell_D\}, \ell_D; \ul{\bL}) & \Cok(j). \ar[l]_-{\cong} &
}\]
Here, the map $\bdry^{alg}$ is a monomorphism induced from the connecting map in the Mayer--Vietoris sequence for $D$-equivariant $\bL$-homology. The map $\bdry^{top}$ is induced from the connecting map in Cappell's exact sequence in $L$-theory. By \cite[Theorem~2]{Cappell_unitary}, the map $c$ is injective. By \cite[Theorem~5(ii)]{Cappell_unitary}, the middle row is exact. By Bartels' theorem \cite{Bartels_splitting}, the bottom horizontal map is an isomorphism. For general group-theoretic reasons, the middle column is exact and $\kappa$ is surjective.

Note, by the calculation in \cite[Theorem~13A.1]{Wall_book}, that $\Ker(k)=0$ for all $n$. Then $c$ is surjective and $\Cok(i)=0$. Therefore $\kappa$, hence $\phi$, is an isomorphism.
\end{proof}

\begin{rem}\label{rem:TSO}
By \fullref{prop:Kflat} and the $s$-cobordism theorem, there is a bijection
\[
\xymatrix{\mcS_\TOP(\ol{X},\bdry\ol{X}) \ar[r]^-{\approx} & \mcS_\TOP^h(\ol{X},\bdry\ol{X}).}
\]
For any compact manifold $W$ of dimension $n \geq 6$, Ranicki  gives a bijection from the geometric structure set to the algebraic structure group \cite[Theorem~18.5]{Ranicki_TOP}:
\[
\xymatrix{
\mcS^h_\TOP(W,\bdry W) \ar[r]^-{\approx} & \mcS^h_{n+1}(W,\bdry W).}
\]
In the case $(W,\bdry W) = (\ol{X},\bdry\ol{X})$ this bijection is valid for $n = 5$ since all the fundamental groups are good in the sense of Freedman--Quinn \cite{FQ}, and it is valid for $n = 4$ since we use homology equivalences on the 3-dimensional boundary.

The Atiyah--Hirzebruch spectral sequence shows that $H_{n+1}(W,\bdry W; \bL/\bL\connective) = 0$ for a compact $n$-manifold $W$ with boundary.   A diagram chase using this fact together with the definitions of the algebraic structure groups give a monomorphism
\[
\xymatrix{\mcS^h_{n+1}(W,\bdry W)~~ \ar@{>->}[r] & \mcS^{per,h}_{n+1}(W,\bdry W).}
\]
\end{rem}

Recall the definitions of $\bdry^{Wall}$ and $\beta$ from the introduction of \fullref{sec:calculation}.

\begin{prop}\label{prop:beta}
The following composite function is a bijection of pointed sets:
\[
\bdry^{Wall} \circ \beta ~:~ \bigoplus_{D \in (\midd)(\Gamma)} \UNil_{n+1}(\Z;\Z^\eps,\Z^\eps) ~\longrightarrow~ \mcS_\TOP(\ol{X},\bdry\ol{X}).
\]
\end{prop}

\begin{proof}
Consider the commutative diagram
\[\small\xymatrix{
\bigoplus \UNil_{n+1}(\Z;\Z^\eps,\Z^\eps) \ar[r] \ar[d]^-{\cong} & \bigoplus L^h_{n+1}(D,w_n) \ar[r] \ar[dl] & L^h_{n+1}(\Gamma,w_n) \ar[r] \ar[d] & \mcS_\TOP(\ol{X},\bdry\ol{X}) \ar[d]^-{\cong} \\
\bigoplus H^D_{n+1}(\{\ell_D\},\ell_D ;\ul{\bL}) \ar[r]^-{\cong}  & H^\Gamma_{n+1}(E,\R^n ;\ul{\bL}) \ar[r]^-{\cong}  & H^\Gamma_{n+1}(\{\R^n\},\R^n; \ul{\bL}) \ar[r]^-{\cong}  & \mcS^{per,h}_{n+1}(\ol{X},\bdry\ol{X}).
}\]
The composition of the three maps in the top row is $\bdry^{Wall} \circ \beta$.
By \fullref{lem:UNil_relhom}, the leftmost vertical map is an isomorphism. 
By \fullref{lem:uncontrolledcalculation}, the leftmost map of the bottom row is an isomorphism.
By \fullref{lem:model} and the Farrell--Jones Conjecture \cite{BL_CAT0}, the middle map of the bottom row is an isomorphism.
By \fullref{lem:split}, the rightmost map of the bottom row is an isomorphism.

In particular, the composite from the upper left of the diagram to the lower right of the diagram must be surjective.
Hence the right vertical map is surjective.
By \fullref{rem:TSO}, the rightmost vertical map is injective.  Hence the right vertical map is bijective.
Therefore, the three top horizontal maps compose to a bijection.
\end{proof}

Observe that the map $\mcS^h_{n+1}(\ol{X},\bdry\ol{X}) \longrightarrow \mcS^{per,h}_{n+1}(\ol{X},\bdry\ol{X})$ from the non-connective structure group to the connective structure group is an isomorphism in our case, since the rightmost vertical map in the diagram above is an isomorphism.

\begin{proof}[Proof of \fullref{thm:calculation}]
This follows from \fullref{prop:alpha}, \fullref{prop:beta}, and \eqref{thm:BCDR}.
\end{proof}

%------------------------------------------------------------------------------
\subsection{Verification of the Topological Rigidity Conjecture}\label{proofTRC}

Lastly, we show that our family of crystallographic examples satisfies our conjecture.

\begin{proof}[Proof of \fullref{toprc} for $(X^n,\Gamma)=(\R^n,\Gamma_n)$]
Consider the commutative diagram
\[\xymatrix{
H^\Gamma_{n+1}(\R^n;\ul{\bL}) \ar[r]^-{A_\fin^\vc} & H^\Gamma_{n+1}(E;\ul{\bL}) \ar[r]^-{A_\vc^\all} \ar[d] & L^h_{n+1}(\Gamma,w_n) \ar[r]^-{\bdry^{Wall}} \ar[d] & \mcS_\TOP(\ol{X},\bdry\ol{X}) \ar[d]^-{\approx}\\
& H^\Gamma_{n+1}(E,\R^n ;\ul{\bL}) \ar[r]^-{\cong} & H^\Gamma_{n+1}(\{\R^n\},\R^n; \ul{\bL}) \ar[r]^-{\cong} & \mcS^{per,h}_{n+1}(\ol{X},\bdry\ol{X}).
}\]
The three bijections hold as in the earlier diagram. Define a precursor map $\widehat{\bdry}$ by
\[
\widehat{\bdry} ~:=~ u^{iso} \circ \Phi \circ \bdry^{Wall} \circ A_\vc^\all ~:~ H^\Gamma_{n+1}(E; \ul{\bL}) \longrightarrow \mcS^{iso}_{rel}(\R^n,\Gamma).
\]
By a theorem of A.~Bartels \cite{Bartels_splitting}, there is a short exact sequence
\[\xymatrix{
0 \ar[r] & H^\Gamma_{n+1}(\R^n;\ul{\bL}) \ar[r]^-{A_\fin^\vc} & H^\Gamma_{n+1}(E;\ul{\bL}) \ar[r] & H^\Gamma_{n+1}(E,\R^n; \ul{\bL}) \ar[r] & 0.
}\]
Then $\widehat{\bdry}$ induces a map $\bdry$ from $\Cok(A_\fin^\vc)$. Using the identification \eqref{iso:Quinn}, we obtain:
\[
\bdry ~:~ H^\Gamma_{n+1}(E_\vc\Gamma; \ul{\bL}/\ul{\bL}_\fin) ~\longrightarrow~ \mcS^{iso}_{rel}(\R^n,\Gamma).
\]
Therefore, since $u^{iso}$ and $\Phi$ are bijections, this desired map $\bdry$ is also a bijection.
\end{proof}

%------------------------------------------------------------------------------
%------------------------------------------------------------------------------
\section{Classification of involutions on tori}\label{sec:classify_involutions}

The goal of this section is to prove \fullref{thm:torus}.

\begin{proof}[Proof of \fullref{thm:torus}\eqref{torus_discrete}]
This is immediate from \fullref{inv_on_torus}(1).
\end{proof}

\begin{proof}[Proof of \fullref{thm:torus}\eqref{torus_rigid}]
The case $n=0$, is trivial: $T^0=\R^0/\Z^0 = pt$.

Assume $n=1$. Set 
$
D^1_\pm ~:=~ \{z=x+iy \in S^1 \subset \mathbb{C} \mid \pm y \geq 0 \}.
$

Write $a,b \in N$ for the fixed points of $\sigma$. Let $f\co D^1_+\to N$ be a homeomorphism of $D^1_+$ onto either  arc in $N$ with endpoints  $a$ and $b$. Extend $f$ to a continuous map $f\co  S^1\to N$ by setting
\[
f(z) = \sigma f (\overline{z})\qquad \forall z\in D^1_-.
\]
Then $f\co  (S^1,C_2) \to (N,C_2)$ is an equivariant homeomorphism.

Assume $n=2,3$.  There is a homeomorphism $f\co  N \to T^n$ (by work of Perelman \cite{Anderson} for $n=3$).    We want to show that each fixed point $x\in N^{C_2}$ has an invariant neighborhood $D$ such that $(D,C_2)$ is homeomorphic to $(D^n,C_2)$, the orthogonal action fixing  only $0$. 

To see this, note the involution $\sigma$ of $(N,x)$ lifts to an involution of the universal cover $(\tN,\tx)$ (for any point $\tx$ over $x$) whose one point compactification provides an involution $\widetilde{\sigma}$ with two fixed points on $S^n$. If this involution is   standard, this yields arbitrarily small standard disk neighborhoods of $\tx$ and the required invariant standard disk neighborhood $(D,C_2)$ of $x$ in $N$. 

But this involution on $S^n$ \emph{is} standard. For, when $n=2$ this was proved by K\'er\'ekjart\`o, Brouwer, and Eilenberg (see \cite{ConstantinKolev_S2});  when $n=3$ it was proved by Hirsch--Smale and  Livesay (see \cite{Rubinstein_S3}).
 
Around each fixed point remove the interior of such an invariant standard disk, to obtain a compact manifold with a free involution, $(N^n_0, \sigma_0)$ whose boundary consists of $2^n$ copies of $S^{n-1}$ with the antipodal involution. This manifold with free involution is smooth. This is by Moise \cite[Thm.~9.1]{Moise}  and Whitehead \cite{Whitehead_smoothing} if $n=3$, and by the classification of surfaces in $n=2$. Gluing back the $2^n$ standard disks, we conclude $N$ is smooth, and $\sigma$ is smooth. 
 
If $n=3$,  a theorem of Meeks--Scott \cite{MS_3manifolds} then proves  there is a flat, invariant Riemannian metric on $(N,C_2)$. So we may assume $N= T^3$ and $C_2$ acts by isometries, and the origin is an isolated fixed point. The group of    all isometries fixing the origin is $O(3)\cap GL_3(\Z)$. Only $-I$ acts with the origin as an isolated fixed point. This is the standard involution on $T^3$.  This proves the theorem when $n=3$.
 
If $n=2$,  we see by the Euler characteristic that $N/C_2$ must be $S^2$, and $(N, C_2)$ must be the two-fold cover, branched at four points of $S^2$. This, again, is the standard involution on $T^2$. This proves the theorem when $n=2$.
 
Assume $n \geq 4$ and $n \equiv 0,1 \pmod{4}$.  By \fullref{inv_on_torus}(2), there is a $C_2$-homotopy equivalence $J\co  N \to T^n$. Recall, from \fullref{isostruc}, the bijection $u\co  \mcS_\TOP(T^n,C_2) \to \mcS_\TOP(\R^n,\Gamma_n)$. By \fullref{lem:nobad}(2), there is a bijection $\chi\co  \mcS_\TOP(\R^n,\Gamma_n) \to \mcS(\Gamma_n)$. By \fullref{thm:calculation}, $\mcS(\Gamma_n)$ is a singleton. Thus $\mcS_\TOP(T^n,C_2)$ is also. Therefore $J$ is $C_2$-homotopic to a homeomorphism.
\end{proof}
 
The proof of \fullref{thm:torus}\eqref{torus_nonrigid} will take a little preliminary work. 
Let $\mcT_n$ denote the set of equivariant homeomorphism classes of $H_1$-negative  $C_2$-manifolds $(N^n, C_2)$ for which $N^n$ has the homotopy type of $T^n$.  
We must prove $\mcT_n$ is infinite if $n\equiv 2,3 \pmod{4}$ and $n\geq 6$. 
Write $\hAut(T^n, C_2)$ for the group of $C_2$-homotopy  classes of $C_2$-homotopy equivalences, $f\co  (T^n, C_2)\to (T^n, C_2)$.  Note, by \fullref{inv_on_torus}\eqref{iso_homotopy_unique}, that
\begin{equation}\label{T}
\mcT_n ~\approx~ \mcS_\TOP(T^n,C_2) / \hAut(T^n,C_2).
\end{equation}

We begin by constructing a homomorphism $\Aut(\Gamma_n)\to \hAut(T^n, C_2)$.  

Recall $\Gamma_n= \Z^n\rtimes_{-1}C_2$. 
For increased clarity below, we shall write $A_n$ for the subgroup of translations in $\G$ and write $\sigma_0 \in \G$ for the reflection through $0$:  
\[
A_n := \{(x,1)\in \Z^n\rtimes_{-1} C_2  \mid x\in \Z^n\} \qquad \sigma_0 := (0, \sigma)\in \Z^n\rtimes_{-1} C_2.
\]
Also, below we shall use the following identification:
\[
k\co  A_n \longrightarrow H_1(T^n) ~;~ \quad (x,1) \longmapsto ([0,1] \to T^n ~;~ \theta \mapsto \theta x + \Z^n).
\]

For each automorphism $a\co  \Gamma_n  \to \Gamma_n $   choose an $a$-equivariant continuous map  
\begin{equation}\label{aeq}
\tJ_a\co  \R^n \longrightarrow \R^n, \quad\text{that is,}\quad \tJ_a(\gamma\cdot v) = a(\gamma)\cdot \tJ_a(v), \enspace \forall (\gamma, v)\in \Gamma_n\times \R^n.
\end{equation}
Note $\tJ_a$  is unique up to $a$-equivariant homotopy.  Since $a(A_n)=A_n$,   we see  $\tJ_a$ descends to a map, $J_a\co (T^n, C_2)\to (T^n, C_2)$. So $[J_a]\in \hAut(T^n, C_2)$. From \eqref{aeq} we see that for  all $a, b\in \Aut(\Gamma_n), \; \tJ_{ab} $ and $\tJ_a\tJ_b$ are $ab$-equivariantly homotopic. 

For  each  $x\in \Gamma_n$,  we write $c(x)$ for the automorphism:
\[
c(x)\co  \Gamma_n \longrightarrow \Gamma_n ~;~ \quad \gamma \longmapsto  x \gamma x^{-1}.
\]
If $t\in A_n$ is any translation, a valid choice for $\tJ_{c(t)}$ is: 
\[
\tJ_{c(t)}\co  \R^n \longrightarrow \R^n ~;~ \quad  v \longmapsto t\cdot v
\]
since \eqref{aeq} holds for this choice. So this construction specifies a  homomorphism
\[
J\co  \Aut(\Gamma_n)/c(A_n) \longrightarrow \hAut(T^n,C_2) ~;~ \quad  [a] \longmapsto [J_a].
\]
Write $\Aut(\G)_{\sigma_0} :=\{a\in \Aut(\G) \mid a(\sigma_0)= \sigma_0\}$. For   $a\in \Aut(\G)_{\sigma_0}$, let $T\co \R^n\to \R^n$ be the unique linear isomorphism   satisfying: $T(t\cdot 0) = a(t)\cdot 0, \; \forall t\in A_n$.  Then $T(\gamma\cdot x)= a(\gamma)\cdot x$ for all $ x\in \R^n$, and therefore a valid choice for $\tJ_a$ is: $\tJ_a=T$.

\begin{prop}\label{kiso}
$J$ is  an isomorphism.
\end{prop}

\begin{proof}
Let $[a]\in \Ker(J)$. We show $[a]=1$.
Using the isomorphism $k\co  A_n\cong H_1(T^n)$
and the fact that $(J_a)_*= \id\co  H_1(T^n)\to H_1(T^n)$, note that $a(t) = t$ for all $t\in A_n$.
Also $J_a$ fixes the discrete set $(T^n)^{C_2}$, since $[J_a]= 1$. Therefore $\tJ_a(0)\in \Z^n$.  Replacing $a\in [a]$ with $a\cdot c(t)$ for a suitable $t\in A_n$, we conclude for our new $a$ that $\tJ_a(0)=0$. So $a(\sigma_0)= \sigma_0$. But $\Gamma_n= \langle A_n, \sigma_0\rangle$. So $a=\id_{\Gamma_n}$ and $J$ is injective.

Now we show $J$ is surjective.
Let $[f]\in \hAut(T^n, C_2)$. Here   $f\co  T^n \to T^n$ is a  $C_2$-map. Let $a''\in \Aut(\G)_{\sigma_0}$ satisfy: $a''(t)= k^{-1}f_*^{-1}(k(t))\in A_n$ for all $t\in A_n$. Then $(J_{a''}f)_*= \id_{H_1(T^n)}$. Let $\tf\co  \R^n\to \R^n$ be a lift of $f$. Note $0\in \R^n_{sing}$, and so  $p:= \tJ_{a''}\tf(0) \in \R^n_{sing}$. There is an involution $\sigma_p\in \G-A_n$ fixing $p$.
   
Let $\res\co  \Aut(\Gamma_n)\to \Aut(A_n)$ be  the restriction homomorphism. Observe $\Ker(\res)$ acts transitively on $\Gamma_n-A_n$, since $\Gamma_n = \langle A_n, \sigma\rangle$ for any $\sigma\in \Gamma_n-A_n$, and each such $\sigma$ is an involution. Therefore  there exists $a'\in \Ker(\res)$ such that $a'(\sigma_p)= \sigma_0$. So $0= \tJ_{a'}(p)=\tJ_{a'}\tJ_{a''}\tf(0)$. Note $\tJ_{a'a''}\tf\co (\R^n, 0)\to (\R^n,0)$ is $A_n$ equivariant. But $J_{a'a''}f$ is $C_2$-equivariant, so $\tJ_{a'a''}\tf$ is $\sigma_0$ equivariant. Therefore $\tJ_{a'a''}\tf$ is $\G$-equivariant and so $\tJ_{a'a''}\tf$ is $\G$-homotopic to $\id_{\R^n}$. Therefore $[f]= [J_a]$, where $a= (a' a'')^{-1}$.
So  $J$ is surjective.
\end{proof}

\begin{proof}[Proof of \fullref{thm:torus}\eqref{torus_nonrigid}]
Assume $n\equiv 2,3 \pmod{4}$ and $n\geq 6$. 
For any group $G$, we are going to abbreviate
\[
H^G ~:=~ H^{G}_{n+1}(E_{\vc}G; \ul{\bL}/\ul{\bL}_{\fin}).
\]
From \fullref{kiso} and \fullref{proofTRC} and \eqref{T}, we see that $H^{\Gamma_n}/\Aut(\Gamma_n) \approx \mcT_n$. So we must prove that this set $H^{\Gamma_n}/\Aut(\Gamma_n)$ is infinite.  
The proof is based on the fact that, for any maximal infinite dihedral subgroup $D$, we have $H^D\cong \UNil_{n+1}(\Z; \Z^\eps, \Z^\eps)$ by \fullref{lem:UNil_relhom}, and so $H^D$ is an infinite set by \eqref{thm:BCDR}.
We will produce an injective map from an infinite set:
\[
H^D/\Aut(D) \longrightarrow H^{\G}/\Aut(\G).
\]
  and this will show that $\mcT_n$ is an infinite set.
  
Take  $D\in (\midd)(\G)$, a maximal dihedral subgroup.
The inclusion $i_D\co  D\to \Gamma_n$  induces a map, $i_{D*}\co  H^D \to H^{\Gamma_n}$. 
It is easy to see that each automorphism of $D$ extends to an automorphism of $\Gamma_n$. 
For  any $a\in \Aut(\Gamma_n)$, we have
\[
a_*\circ i_{D*} ~=~ i_{a(D)*}\circ(a|_{D})_* ~:~ H^D \longrightarrow H^{\Gamma_n},
\]
where $a|_D\co  D\to a(D)$ denotes the restriction of $a$.

By \fullref{lem:uncontrolledcalculation}, the induced map $i_{D*}\co H^D\to H^\G$ is a monomorphism, and $i_{D*}(H^D) \cap i_{a(D)*}(H^D)=0$ if $a(D)$ is not conjugate to $D$. But if $a(D)=c_\gamma(D)$ for some $\gamma\in \Gamma_n$, then the following diagram commutes:
\begin{equation}\label{diagram}\xymatrix{
H^D \ar[rr]^-{i_{D*}} \ar[d]_-{(c_{\gamma^{-1}}\,  a|_D)_*} & & H^{\Gamma_n} \ar[d]^-{(c_{\gamma^{-1}})_*a_*}\\
H^D \ar[rr]^-{i_{D}*} & & H^{\Gamma_n}.
}\end{equation}

By a theorem of Taylor \cite{Taylor_inner}, for any $\gamma\in D$, the map $c_{\gamma *}\co  H^D\to H^D$ is just multiplication by $(-1)^k$, if $w(\gamma)=(-1)^k$.  Now if $n\equiv 3 \pmod{4}$, then $w(\gamma)=1$ for all $\gamma$, and if $n\equiv 2 \pmod{4}$, then $2\cdot H^D=0$, so, in all cases $c_{\gamma *} = \id$. For the same reason, if $\gamma\in \G$, then $c_{\gamma *} = \id\co  H^{\G}\to H^{\G}$.

This together with \eqref{diagram} proves, first, that  the induced map  
\[
(i_D)_* \co  H^D/\Aut(D)\to H^{\G}/\Aut(\G) 
\]
is injective,  and second that $\Inn(D)$ acts trivially on the infinite group $H^D$. But $\Aut(D)/\Inn(D)\cong C_2$, so $H^D/\Aut(D)$ is an  infinite set.
Therefore   $H^{\Gamma_n}/\Aut(\Gamma_n)$ and $\mcT_n$ are also infinite sets, as required.
\end{proof}

%------------------------------------------------------------------------------
%------------------------------------------------------------------------------
\section{A non-trivial element of $\mcS(\Gamma_n)$}\label{sec:quick}

In this   section, which is independent of the rest of the paper, we give a classical argument for the existence of non-trivial elements of $\mcS (\Gamma_n)$ for some $n$. Indeed the argument for  the case  $n\equiv 2 \pmod 4$,  could have been written in 1976.  It was in fact pointed out by   Weinberger  to the first author many years ago.

\begin{thm}\label{quickc} Suppose $ n\equiv 2 $ or $ 3 \pmod{4}$  and $n\geq 6$. There exists a cocompact action of the group $\Gamma_n$ on a manifold $M^n$ such that $(M^n,\Gamma_n)$ is simply isovariantly homotopy equivalent to $(\R^n, \Gamma_n)$ but is not equivariantly homeomorphic to $(\R^n, \Gamma_n)$. 
\end{thm}
The proof depends only on  an idea of Farrell \cite{Farrell_exponent} and Cappell's Splitting Theorem (\cite[Theorem 6]{Cappell_free}, \cite{Cappell_split}). It does not depend on  \cite{BL_CAT0}.  

Let $w_n\co \Gamma_n = \Z^n \rtimes C_2 \to \{\pm1\}$ be the homomorphism such that $\ker(w_n)= \Z^n$ if $n$ is odd, and $w_n(\G)=\{1\}$ if $n$ is even.
By \fullref{amalg}, there is a group isomorphism $\Gamma_{n-1} *_{\Z^{n-1}}\Gamma_{n-1}\to \Gamma_n$. By Cappell \cite{Cappell_unitary}, this decomposition defines a split monomorphism
\[
\rho\co  \UNil_{n+1}^h(R; \mcB, \mcB') \to L_{n+1}^h(\Z[\Gamma_n], w_n)~;~ R= \Z[\Z^{n-1}] ~,~ \mcB = \mcB' = \Z[\Gamma_{n-1}-\Z^{n-1}].
\]
Here $R$ is a ring with involution given by: $\overline{a}= a^{-1}$ for all $a\in \Z^{n-1} \subset \Gamma_{n-1}$. Also $\mcB$ and $\mcB'$ are $R$-bimodules with involution: $\overline{b} =( -1)^nb^{-1}$ for all $b \in \Gamma_{n-1}-\Z^{n-1}$.   

\begin{lem}[Cappell]\label{cap}
The action of the abelian group $L_{n+1}^h(\Z[\Gamma_n],w_n)$ on the set $\mcS_\TOP(\ol{X},\bdry\ol{X})$ restricts to a free action of $\UNil_{n+1}^h(R;\mcB ,\mcB')$ on $\mcS_\TOP(\ol{X}, \bdry\ol{X})$.
\end{lem}    

\begin{proof}
Cappell's Splitting Theorem (see \cite[Thm.~6]{Cappell_free}, \cite{Cappell_split}) applies only to a \emph{closed} manifold $X$  with $\Gamma_n= \pi_1(X)$, if $X$ admits a splitting $X=X_1\cup_{Y}X_{2}$  consistent with the decomposition $\Gamma_n\cong \Gamma_{n-1}*_{\Z^{n-1}}\Gamma_{n-1}$.  As stated, it does not apply to $\ol{X}$, since $\bdry\ol{X}$ is nonempty.  But  $\ol{X}$ does split along a \emph{closed submanifold}:
\[
\ol{X}= X_1\cup_Y X_2\quad\text{where}\quad Y:=\left\{[t_1, \dots, t_n]\in T^n \big| \; t_1 =\pm  {\textstyle\frac{1}{4}} \right\}/C_2 \,.
\]
Here $X_1$ (and $X_2$) is defined similarly in $\ol{X}$ but with $t_1 \in [-\frac{1}{4},\frac{1}{4}]$ (respectively, $t_1 \in [\frac{1}{4},\frac{3}{4}]$).
The fundamental groups of $X_1, X_2, Y$ are the groups  appearing in \fullref{amalg} with $f(a_1, \dots , a_n) = a_1$,
So, by \fullref{cap_app},  we are done. Again, the key point is $\bdry Y=\varnothing$.
\end{proof}

\begin{lem}\label{unonzero}
$\UNil_n(\Z; \Z^{\eps},\Z^{\eps})$ is a   summand of $\UNil_n(R; \mcB, \mcB')$, if $\eps=(-1)^n$. Furthermore, the groups
$\UNil_{4k+2}(\Z; \Z, \Z)$ and $\UNil_{4k+3}(\Z; \Z,\Z)$ are non-zero.
\end{lem}

\begin{proof} 
The first claim is immediate from the split epimorphism $$\varepsilon \co  \G \to \Gamma_1 = C_2 * C_2$$ of \fullref{amalg}, which induces a split epimorphism 
$$
\UNil_n(R; \mcB, \mcB') \to \UNil_n(\Z; \Z^{\eps},\Z^{\eps}).
$$

The second claim has been known for many years.  See \cite{Cappell_splittingobst} or \cite{CK_nilgps}. But for the reader's convenience, here is a
very easy proof. 
Farrell (see \cite{Farrell_exponent}) extended Cappell's homomorphism $\rho$, mentioned above, to a homomorphism,
\[
\rho'\co  \UNil_{2k}(R; R, R)\to L_{2k}(R)
\]
for any ring with involution $R$. But  the  non-zero element of $L_2(\Z)$ is the class of the rank-two $(-1)$-quadratic form with Arf invariant 1. This element is  $\rho'([\zeta])$
where $[\zeta]\in \UNil_2(\Z;\Z,\Z)$ is the class of the unilform, $\zeta=(P,\lambda, \mu, P', \lambda', \mu')$, where,
\[
P=\Z e,\; P'=\Z f,\;  \lambda = 0,\; \lambda'=0; \quad \mu(e)=\mu '(f)=1\pmod{2}.
\]
Finally a quick proof that $\UNil_{4k+3}(\Z; \Z,\Z)$ has an element of order $4$ can be found in \cite[Corollary~1.9]{CD}. It uses almost no machinery.
\end{proof}

\begin{proof}[Proof of \fullref{quickc}]
By \fullref{unonzero}, \fullref{cap}, and \fullref{lem:compactification}, there is an element $[M,f]\neq [\R^n, \id]$ in $ \mcS_\TOP^{iso}(\R^n,\Gamma_n)$. By \fullref{lem:nobad}, we conclude that $(M, \G)$ is not equivariantly homeomorphic  to $(\R^n,\Gamma_n)$.
\end{proof} 

%------------------------------------------------------------------------------
%------------------------------------------------------------------------------
\subsection{Free action of $\UNil$ on the structure set of a pair} \label{unil_appendix}
\

Our purpose here is to show that our \fullref{cap} is a formal consequence of the  $L$-theoretic exact sequence of Cappell, appearing in \cite{Cappell_unitary, Cappell_split}.

Let $X$ be a compact, connected topological manifold of dimension $n \geq 6$. Let $Y$ a connected, separating, codimension-one submanifold of $X$ without boundary (that is, $\bdry Y = \varnothing$). Assume the induced map $\pi_1(Y) \to \pi_1(X)$ of groups is injective. Furthermore, assume the induced map $\pi_1(\bdry X) \to \pi_1(X)$ of groupoids is injective.

Write $X = X_1 \cup_Y X_2$ for the induced decomposition of manifolds. Write $G = G_1 *_F G_2$ for the induced injective amalgam of fundamental groups, where $F := \pi_1(Y)$. Finally, write $H := \pi_1(\bdry X)$ as the fundamental groupoid of the boundary.

For simplicity of notation, we shall suppress all the orientation characters. Furthermore, to avoid $K$-theoretic difficulties, \emph{we assume throughout this subsection} that the projective class group for the codimension-one submanifold $Y$ vanishes:
\[
\tK_0(\Z[F]) ~=~ 0.
\]

\begin{thm}\label{cap_app}
On the structure set $\mcS^h_\TOP(X,\bdry X)$ of the pair, Wall's action of the group $L_{n+1}^h(G)$  restricts to a free action of Cappell's subgroup
\[
\UNil_{n+1}^h ~:=~ \UNil_{n+1}^h(\Z[F]; \Z[G_1-F],\Z[G_2-F]).
\]
\end{thm}

Thus we slightly generalize the case of $\bdry X=\varnothing$ of Cappell \cite[Thm.~2]{Cappell_free}.
Our proof relies only on his algebraic results \cite[Thm.~2, Thm.~5]{Cappell_unitary} \cite{Cappell_split}.

\begin{thm}[Cappell]\label{thm:MVCappell}
There is a homomorphism
\[
\iota ~:~ \UNil^h_* ~\longra~ L^h_*(G)
\]
whose composite with a map of Wall (see \cite[Thm.~9.6]{Wall_book}) is an isomorphism:
\[
\UNil^h_* ~\xrightarrow{~\iota~}~ L^h_*(G) ~\longra~ L^h_*
\left(\begin{gathered}\xymatrix{
F \ar[r] \ar[d] & G_1 \ar[d]\\
G_2 \ar[r] & G
}\end{gathered}\right).
\]
Furthermore, there is an exact sequence
\[
\cdots \to \UNil^h_{*+1} \oplus L^h_*(F) \to L^h_*(G_1) \oplus L^h_*(G_2) \to L^h_*(G) \xrightarrow{\left(\begin{smallmatrix}s\\ \bdry\end{smallmatrix}\right)} \UNil^h_* \oplus L^h_{*-1}(F) \to \cdots
\]
such that $s \circ \iota = \id$. In particular, $\iota$ is split injective with a preferred left-inverse.
\qed
\end{thm}

Using the first part of \fullref{thm:MVCappell}, the first two authors formally identify $\UNil$ with relative equivariant homology \cite{CD2}. Write $\all$ as the family of all subgroups of $G$. Write $\fac$ as the family of subgroups of $G$ conjugate into either $G_1$ or $G_2$.

\begin{thm}[Connolly--Davis]\label{thm:CD}
The following composite map is an isomorphism:
\[
\phi ~:~ \xymatrix{\UNil^h_* \ar[r]^-{\iota} & L^h_*(G) ~=~ H^G_*(E_\all G; \ul{\bL}^h) \ar[r]^-{h} & H^G_*(E_\all G, E_\fac G; \ul{\bL}^h).}
\]
\end{thm}

This implies a relative version; recall $H = \pi_1(\bdry X)$. Consider the homomorphism
\[
j ~:~ L^h_*(G) ~\longra~ L^h_*(G,H).
\]
Write the fundamental groupoid $H = H_1 \sqcup \cdots  \sqcup H_m$ as the disjoint union of its vertex groups $H_i$. The associated $G$-set is defined by
\[
G/H ~:=~ G/H_1 \sqcup \cdots \sqcup G/H_m.
\]
Observe, since $Y \cap \bdry X = \varnothing$, for each $i$ that $H_i \subset G_1$ or $H_i \subset G_2$; hence $H_i \in \fac$.  Therefore there is a canonical $G$-map $G/H \to E_\fac G$.

\begin{cor}\label{cor:relUNil}
There is a split short exact sequence
\[\xymatrix{
0 \ar[r] & H^G_*(E_\fac G, G/H; \ul{\bL}^h) \ar[r]^-{A_\fac} & L^h_*(G,H) \ar[r]^-{s'} & \UNil^h_* \ar[r] & 0.
}\]
The preferred right-inverse for $s'$ is the composite $j \circ \iota$.
\end{cor}

\begin{proof}
There is a long exact sequence of the triple:
\[
\cdots \to H^G_*(E_\fac G, G/H; \ul{\bL}^h) \to H^G_*(E_\all G, G/H; \ul{\bL}^h) \xrightarrow{k} H^G_*(E_\all G, E_\fac G; \ul{\bL}^h) \to \cdots.
\]
So, by \fullref{thm:CD}, we may define a homomorphism
\[
s' ~:~ 
L^h_*(G,H) ~=~ \xymatrix{H^G_*(E_\all G, G/H; \ul{\bL}^h) \ar[r]^-{k} & H_*(E_\all G, E_\fac G; \ul{\bL}^h) \ar[r]^-{\phi^{-1}} & \UNil^h_*.}
\]
That is, $s' := \phi^{-1} \circ k$. Note $h = k \circ j$. Recall $\phi = h \circ \iota$.
Then
\[
s' \circ (j \circ \iota) ~=~ (\phi^{-1} \circ k) \circ (j \circ \iota) ~=~ \phi^{-1} \circ (h \circ \iota) = \phi^{-1} \circ \phi = \id.
\]
Therefore, $k$ has right-inverse $j \circ \iota \circ \phi^{-1}$, and the above exact sequence splits.
\end{proof}

Now we are ready to prove the main theorem of this subsection.

\begin{proof}[Proof of \fullref{cap_app}]
Ranicki defined algebraic structure groups $\mcS^h_*(X, \bdry X)$, a homomorphism $L^h_*(G,H) \to \mcS^h_*(X,\bdry X)$, and a pointed bijection
\[
\xymatrix{\mcS^h_\TOP(X,\bdry X) \ar[r]^-{\approx} & \mcS^h_{n+1}(X,\bdry X)}
\]
such that it is equivariant with respect to the actions of $L^h_{n+1}(G,H)$ (see \cite{Ranicki_TOP}).
Also observe, from \fullref{rem:TSO}, that there is a monomorphism
\[
\xymatrix{\mcS^h_{n+1}(X,\bdry X)~~ \ar@{>->}[r] & \mcS^{per,h}_{n+1}(X,\bdry X).}
\]
Write $W$ as the composite homomorphism $L^h_{n+1}(G,H) \to \mcS^{per,h}_{n+1}(X,\bdry X)$, which is compatible with Wall's action of $L^h_{n+1}(G,H)$ on the structure set $\mcS^h_\TOP(X,\bdry X)$. Thus it suffices to show that the following composite is a monomorphism:
\[\xymatrix{
\UNil^h_{n+1} \ar[r]^-{\iota} & L^h_{n+1}(G) \ar[r]^-{j} & L^h_{n+1}(G,H) \ar[r]^-{W} & \mcS^{per,h}_{n+1}(X,\bdry X).
}\]

By definition of the algebraic structure groups, there is an exact sequence
\[\xymatrix{
H_{n+1}(X,\bdry X; \bL) \ar[r]^-{A} & L^h_{n+1}(G,H) \ar[r]^-{W} & \mcS^{per,h}_{n+1}(X,\bdry X).
}\]
Also, using \fullref{QR=DL}, there is a commutative diagram of assembly maps:
\[\xymatrix{
H^G_{n+1}(\tX, G \times_H \widetilde{\bdry X}; \ul{\bL}^h) \ar[d]_-{=} \ar[r] & H^G_{n+1}(E_\fac G, G/H; \ul{\bL}^h) \ar[d]_-{A_\fac}\\
H_{n+1}(X,\bdry X; \bL) \ar[r]^-{A} & L^h_{n+1}(G,H).
}\]
Then, by \fullref{cor:relUNil}, note:
\[
\Ker(W) = \Img(A) \subseteq \Img(A_\fac) = \Ker(s') \quad\text{and}\quad  \Img(j \circ \iota) \cap \Ker(s') = 0.
\]
So $W \circ j \circ \iota$ is a monomorphism. Therefore $\UNil^h_{n+1}$ acts freely on $\mcS^h_\TOP(X,\bdry X)$.
\end{proof}
  
%------------------------------------------------------------------------------
%------------------------------------------------------------------------------
\appendix
\section{From equivariance to isovariance}\label{eq_to_iso}
 
We  want to prove that the forgetful map $\psi\co \mcS^{iso}_\TOP(T^n, C_2) \to \mcS_\TOP(T^n, C_2)$ is bijective when $n\geq 4$.   
It seems best to approach this from a   general study  of isovariance. It shall be immediate from \fullref{equi-iso} below that $\psi$ is bijective.   The assumption of a discrete singular set in our \fullref{equi-iso} is key.

Let $G$ be a finite group. For any $G$-spaces $X$ and $Y$, write $[X,Y]_G$ and $[X,Y]_G^{iso}$ for the set of $G$-equivariant and $G$-isovariant homotopy classes of maps, respectively.
 
Let $X$ be a  $G$-space with a fixed point $p$.
The \emph{homotopy link} of $p$ in   $X$, denoted $t^pX$, and the \emph{homotopy tangent space} of $X$ at $p$, denoted $t_pX$,  are defined by:
\begin{eqnarray*} 
t^pX &:=& \Holink(X,p) = \left\{\sigma\co [0,1]\to X \mid \sigma^{-1}(p)= \{0\} \right\}\\
t_pX &:=& t^pX\cup\{\sigma_p\}.
\end{eqnarray*}
Here $\sigma_p$ denotes the constant path at $p$, and
$t_pX$ has the compact-open topology. This is the metric topology of uniform convergence if $X$ is a metric space.

There is a $G$-subspace $X_{(p)} \subset X$ and isovariant evaluation map $e_1$ defined by
\[
X_{(p)} := (X-X^G)\cup\{p\}, \qquad
e_1: t_pX_{(p)} \longra X ~;~ \sigma \longmapsto \sigma(1).
\]

Let $U$ be a neighborhood of  $p$ in $X$.  If  $U$ is   homeomorphic to $\R^n$, then
\begin{align}\label{con}
e_1 \text{ restricts to a homotopy equivalence, } t^pU\simeq \R^n-\{0\}.
\end{align}
If $U$ is $G$-invariant and $p$ is an isolated fixed point, then the inclusion $\iota\co  t_pU\to t_pX$ is an isovariant homotopy equivalence. 

An action of a group $G$ on a set $X$ is {\em semifree} if the action is free away from the fixed set, that is, the action of $G$ on $X-X^G$ is a free action.

\begin{lem}\label{metric}
Let $X$ and $Y$ be metric spaces on which $G$ acts semifreely and isometrically. Assume $q$ is an isolated fixed point of $Y$. The rule $f\mapsto f|_{X-X^G}$ gives a bijection between isovariant and  equivariant homotopy classes, 
\[
[X, t_qY]_G^{iso} \; \cong \; [X-X^G, t^qY]_{G}.
\]
\end{lem}

\begin{proof}
Let $f\co  X-X^G \to t^qY$ be a $G$-map. We  first show $f$ is $G$-homotopic to an extendible map, $f'$. (Here $f'$ is \emph{extendible}, if $\lim_{x\to x_0} f'(x) = \sigma_q$ for any  $x_0\in X^G$.)
This will prove that the restriction map is a surjection. We assume $X$ and $Y$ have metrics, $d_X$ and $d_Y$, bounded by $1$. Write $d_t$ for  the induced metric on  $t_qY$.

For $x\in X-X^G$, set $\|x\| := d_X(x, X^G)$. For $\sigma\in t^qY$, set $\|\sigma\| := d_t(\sigma,\sigma_q)$. If, for all $x\in X-X^G,\; \|f(x)\|\leq \|x\|$, then $f$ is obviously extendible. In $(X-X^G)\times I$, the subset $(X-X^G)\times 0$ is disjoint from the closed subset
\[
B ~:=~ \left\{(x, t)\in (X-X^G)\times I \mid d_Y(f(x)(t), q)\geq \|x\| \right\}.
\]
Consider the continuous map
\[
\phi ~:~ X-X^G \longrightarrow (0,1] ~;~ x \longmapsto d_\times((x, 0), B),
\]
where $d_\times$ denotes the product metric on $(X-X^G)\times I$. 
Observe
\[
\phi(x)= d_\times((x,0), (x, \phi(x))).
\]
If $0\leq t<\phi(x)$ then $(x, t)\notin B$. Therefore $d_Y(f(x)(t), q)\leq \|x\|$ for all $t\in [0, \phi(x)]$.

Define $f'\co X-X^G\to t^qY$ as the map whose adjoint is
\[
\mcA f' := \mcA \circ \Phi \co  (X-X^G)\times I \longrightarrow Y \quad\text{where}\quad \Phi(x, t) := (x, \;t\cdot \phi(x)).
\]
By construction, $\|f'(x)\|\leq \|x\|$ for all $x\in X-X^G$. So $ f'$ is extendible.
But a $G$-homotopy, $f_s, \; 0\leq s\leq 1 $, from $f$ to $f'$ is defined by
\[
\mcA f_s := \mcA f\circ \Phi_s \co  (X-X^G)\times I \longrightarrow Y \quad\text{where}\quad \Phi_s(x, t) := (x, \;t\cdot (s\phi(x)+ (1-s))).
\]

Note that if $f$ is extendible, then each  $f_s$ is extendible too.

The same simple argument shows that if $f\co (X-X^G)\times I\to t^qY$,  is a homotopy between two extendible $G$-homotopy equivalences, then $f'$ supplies an extendible homotopy.  One merely changes the definition of $B$ to
\[
B ~:=~ \left\{ (x, t)\in (X-X^G)\times I \mid d_Y(f(x,s)(t), q)\geq \|x\| \text{ for some } s\in [0,1] \right\}.
\]
This proves that $[X-X^G, t^qY]_G \; \cong \; [X, t_qY]_G^{iso}$, as required.
\end{proof}

Because our argument employs a somewhat unusual form of Poincar\'e duality, namely \eqref{pd}, we spend some space introducing it here.  

\begin{lem}\label{lem:PD}
Let $(X,A)$ be a compact Hausdorff $G$-pair satisfying
\begin{multline}\label{hyp:G_NDR}
\exists \text{ $G$-homotopy $h \co  X \times I \to X$ and  $ \exists ~G$-neighborhood $N$ of $A$ so that } \\
\forall x \in X ~:~ h(x,0) = x, \quad \forall a \in A ~:~ h(a,t) = a, \quad \forall n \in N ~:~ h(n,1) \in A.
\end{multline}
Assume that $X-A$ is an $m$-dimensional manifold with boundary, and that the action of $G$ on $X-A$ is free.

Then for any left $\Z G$-module $B$, there is an isomorphism
\begin{equation}\label{pd}
H^*_G(X-A, \p(X-A); B) ~\cong~ H^G_{m-*}(X, A;B).
\end{equation}
\end{lem}

\begin{note}
For any $G$-pair $(X,A)$, the groups $H^G_*(X,A;B)$ and $ H_G^*(X,A;B)$ denote the homology of the complexes $C(X,A)\otimes_{\Z G}B$ and $\mathrm{Hom}_{\Z G}(C(X,A), B)$. 
\end{note}

\begin{proof}[Proof of \fullref{lem:PD}]
Let $h$ and $N$ be as in \eqref{hyp:G_NDR}.
Let $f \co  X \to X$ be the function $x \mapsto h(x,1)$.
Note $f(X) = A$.
Consider the collection
\[
\mcU ~:=~ \{ G\text{-open } U \text{ in } X \mid A\subset U\subset N \}. 
\]
For each $U\in \mcU$, write $f_U\co  C(X,U)\to C(X,A)$ for the chain map induced by $f$ and $i_U\co  C(X,A)\to C(X,U)$ for the chain map induced by the inclusion $A\subset U$.
 
Since $X-A$ is a non-compact $m$-manifold with boundary, Poincar\'e duality takes the form of an isomorphism
\[
\mcP ~:~ H^*_G(X-A, \p(X-A); B) \longrightarrow H_{m-*}^{G,lf}(X-A;B),
\]
where the right-hand group is the homology of the complex $C^{lf}(X-A)\otimes_{\Z G}B$.
We identify the complex $C^{lf}(X-A)$ of locally finite chains with the inverse limit
\[
C^{lf}(X-A) ~=~ \lim_{U\in\mcU} C(X, U).
\]
The systems of maps given by $f_U$ and $i_U$ then define maps:
\[
\xymatrix{C(X,A) \ar[r]^-{i_\mcU} & \displaystyle\lim_{U\in\mcU} C(X, U)=C^{lf}(X-A) \ar[r]^-{f_\mcU} & C(X,A).}
\]
We will prove below that $i_\mcU$ and $f_\mcU$ are chain homotopy inverse over $\Z G$.
Therefore the required isomorphism is the composite
\[
(i_\mcU\otimes \id_B)_*^{-1}\circ \mcP ~:~ H^*_G(X-A, \p(X-A);B) \longrightarrow H_{m-*}^G(X,A; B).
\]
 
We now show $i_\mcU$ and $f_\mcU$ are chain homotopy inverse to one another. First note that the map $f_\mcU \circ i_\mcU\co C(X,A)\to C(X,A)$ equals the chain map $f_\#$ induced by $f$. So the homotopy $h$ induces a $\Z G$-chain homotopy:  $\id_{C(X,A)} \simeq_G  (f_\mcU\circ i_\mcU)$.
 
Finally, we show that there exists a $\Z G$-chain homotopy $\id_{C^{lf}(X-A)}\simeq_G (i_\mcU\circ f_\mcU)$. 
Choose a sequence $(n \mapsto O_n)$ of open, $G$-invariant neighborhoods of $A$ in $X$ such that $\Cl_X(O_{n+1}) \subset O_n \subset N$ and $A = \bigcap_{n=1}^\infty O_n$. Define an order preserving function, $T\co \mcU\to \mcU$ by setting $T(U)= O_n$ where $n$ is the first integer for which $h(\Cl_X(O_n) \times I) \subset U$.
Write $j_U\co  C(X, T(U))\to C(X,U)$ for the map induced by the inclusion $T(U)\subset U$.
Then for each $U$, the homotopy $h$ induces a chain homotopy
\[
h_U ~:~ j_U \simeq_G (i_U\circ f_{T(U)}) ~:~ C(X,T(U)) \longrightarrow C(X,U).  
\]
But note that:
\[
\lim_{U\in \mcU} j_U = \id_{C^{lf}(X-A)} \qquad \lim_{U\in \mcU} i_U = i_\mcU \qquad \lim_{U\in \mcU} f_{T(U)} = f_\mcU.
\]
Therefore $\lim_{U\in \mcU} h_U $ provides the desired chain homotopy.
\end{proof}
 
Recall that a map $f\co (X,A)\to (Y,B)$ is \emph{strict}  if $f^{-1}(B)= A$.
A homotopy equivalence between the pairs $(X,A)$ and $(Y,B)$ is \emph{strict} if there are corresponding homotopies $(X\times I, A\times I) \to (Y,B)$ and $(Y\times I, B\times I)\to (X,A)$ that are strict.

\begin{defn}
Let $X$ be a $G$-space, and let $A$ a closed $G$-subset of $X$.
We call a pair $(U,h)$ a \emph{tamed neighborhood} if it consists of a $G$-invariant neighborhood $U$ of $A$ in $X$ and a strict $G$-map
$h\co  (U\times I, \; U\times\{0\}\cup A\times I) \to (X, A)$ that restricts to the inclusion on  $U\times\{1\}$ and the projection on $A \times I$.
We say $A$ is \emph{tame in $X$}
if $A$ has a tamed neighborhood $(U, h)$ in $X$.
It follows that $(X,A)$ satisfies \eqref{hyp:G_NDR}.
\end{defn}

\begin{rem}\label{rem:Quinn_HSS}
Note that if $A=\{p\}$ is an isolated fixed point with tamed neighborhood $(U,h)$, then $U^G=\{p\}$
and the adjoint of $h\co U\times I\to X$ is an isovariant map $\lambda_p\co U \to t_pX$.
Furthermore, it follows from Quinn \cite[Prop.~3.6 or Prop.~2.6]{Quinn_HSS} that if $X$ is assumed to be a $G$-manifold and $V$ is any neighborhood of $p$, then there exists a tamed neighborhood $(U,h)$ of $p$ with $U\subset V$. 
\end{rem}

\begin{lem}\label{obstr}
Let $U^m$ be a compact semifree $G$-manifold with $G$-collared boundary. Assume $U^G$ is tame in $U$. Let $V^n$ be a $G$-manifold, and $q$ an isolated fixed point of $V$. Let $f\co U\to t_qV$ be a $G$-map such that $f|_{\p U}$ is isovariant. If $(U, U^G)$ is $1$-connected and $m\leq n+1$, then $f$ is $G$-homotopic $\rel \p U$ to an isovariant map.
\end{lem}

\begin{proof}
There exists a $G$-map $F\co  U-U^G\to t^qV$ extending the $G$-map $f|_{\p (U-U^G)}\co \p (U-U^G)\to t^qV$, since, using \fullref{lem:PD}, the obstructions lie in the groups
\[
H^i_G\left(U-U^G,\p (U-U^G); \pi_{i-1}(t^q V)\right) ~\cong~ H_{m-i}^G\left(U,U^G; \pi_{i-1}(t^q V)\right).
\]
If $i<n$ then the coefficient group is zero, by \eqref{con}.
If $i\geq n$ then the homology group is zero, because $(U, U^G)$ is a $1$-connected pair and $m-i\leq m-n\leq 1$.
Therefore, by \fullref{metric}, there is an isovariant map $f'\co  U\to t_qV$ such that $f'|_{\p U}$ is isovariantly homotopic to $f|_{\p U}$.
But since $\p U$ has a $G$-collar in $U$, the isovariant $G$-homotopy extension property applies, and we can choose $f'$ so that $f'|_{\p U}= f|_{\p U}$.
Finally since $t_qV$ is $G$-contractible, the maps $f, f'\co U\to t_qV$ are $G$-homotopic.
\end{proof}

Let $X$ be a $G$-manifold with boundary.  Recall the \emph{singular set} of $X$ is
\[
X_{sing} ~:=~ \{x\in X \mid gx=x \text{ for some } g \neq 1 \in G \}.
\]
A neighborhood $U$ of $X_{sing}$ in $X$ is a \emph{$k$-neighborhood} if it is a $G$-invariant codimension zero submanifold  with bicollared frontier in $X$, such  that  the pair $(U, X_{sing})$ is $k$-connected.
We only use this concept for $k=0,1$.
 
\begin{thm}\label{equi-iso}
Let $X^n$ and $Y^n$ be compact  $G$-manifolds without boundary.
Assume $X_{sing}$ and $Y_{sing}$ are finite sets, and assume $n\geq 4$.
Let $f\co X\to Y$ be a $G$-map such that the restriction $f_{sing}\co  X_{sing}\to Y_{sing}$ is bijective.    
\begin{enumerate}   
\item
If $f$ is $1$-connected, then $f$ is $G$-homotopic to an isovariant map.
   
\item
Suppose $f=f_0$ is isovariant and $2$-connected.
If $f_0$ is $G$-homotopic to an isovariant map $f_1\co  X\to Y$, then $f_0$ is  $G$-isovariantly homotopic to $f_1$.
\end{enumerate}
\end{thm}

\begin{proof}[Proof of \fullref{equi-iso}(1)]
We argue in two major steps.

\emph{Step 1:} We find $0$-neighborhoods $U, V$ of $X_{sing}, Y_{sing}$ such that $U=f^{-1}(V)$.

It follows from Siebenmann's thesis \cite{Siebenmann_thesis} that each neighborhood of $X_{sing}$ contains a $0$-neighborhood if $n= 4$, and a $1$-neighborhood if $n\geq 5$.
By \fullref{rem:Quinn_HSS}, we can choose a $0$-neighborhood $V$ of $Y_{sing}$ in $Y$, so small that for each $q\in Y_{sing}$, the component $V_q$ containing $q$, admits the structure of a tamed neighborhood of $q$, say $(V_q, h_q)$.
By \cite{FQ} applied to $(f-f_{sing})/G$, we may assume $f$ is transverse to $\p V$.
Then $N:=f^{-1}(\p V)$ is a bicollared codimension one $G$-submanifold of $X$.
It is the boundary and frontier of $U:= f^{-1}(V)$, a $G$-neighborhood of $f^{-1}(Y_{sing})$.

Our desire is that $U$ be a  $0$-neighborhood. We plan to accomplish this by handle exchanges along $N$ realized through a homotopy of $f$.

Define the closures $X_0 := \Cl_X(X-U)$ and $ Y_0 := \Cl_Y(Y-V)$.
Note $X_0$ is a manifold with boundary $N= \p X_0$.
Also note $X=U\cup_N X_0$ and $Y=V\cup_{\p V} Y_0$.

We now recall the aforementioned notion of \emph{handle exchange} along  $N$.
 
Suppose  we can find a map, $i\co  (D^k, \p D^k)\times \{0\}\to (X_0, N)$  (or alternatively, a map, $i\co  (D^k, \p D^k)\times \{0\}\to (U, N)$), together with an extension of $f\circ i$ to a map
\[
j ~:~ (D^k,\p D^k)\times (I,  \{1\}) \longrightarrow (Y_0, \p Y_0) \qquad (\text{or to } (V,\p V) ).
\]
If $k< n/2$, we can, after a homotopy, thicken $i$ to an equivariant embedding and an equivariant extension still called $i$; also, we can thicken $j$ to a continuous $G$-map, still called $j$: 
\begin{align*}
i ~:~ G\times (D^k, \p D^k)\times D^{n-k} \longrightarrow (X_0, N) &\qquad (\text{or to }(U-U^G,N)).\\
j ~:~ G\times (D^k, \p D^k)\times D^{n-k}\times (I,\{1\}) \longrightarrow (Y_0,  N)  &\qquad (\text{or to }(V-V^G,\p V)).
\end{align*}

Now deform $f$ by a $G$-homotopy, stationary off   $  i(G\times \Int(D^k\times D^{m-k}))$,   to a map $f'$ so that $f'$ is still transverse to $N$, but 
\begin{eqnarray*}
f'^{-1}(V) &=& U\cup   i(G\times D^k\times {\textstyle \frac{1}{2}}D^{m-k})\\
(\text{or } f'^{-1}(Y_0) &=& X_0\cup   i(G\times D^k\times {\textstyle \frac{1}{2}}D^{m-k})).
\end{eqnarray*}
Note this homotopy is $\rel \p X$.
If $q\in Y_{sing}$, set $U_q:=f^{-1}(V_q)$ and $N_q:= N\cap U_q$.

Assume $k=1$. This handle exchange process decreases the number of components of $N$ provided that $i$ is chosen so that $\mathrm{Image}(i)/G$ meets two components of $N/G$. After finitely many such handle exchanges  then, we arrive at a map $f'$ for which $N_q$ is connected for each $q\in Y_{sing}$. Therefore  $U_q$ is connected too. So $U$ is a  $0$-neighborhood of $X^G$.

\emph{Step 2:} In this step we find a $G$-homotopy from $f$ to a $G$-isovariant map.

From Step 1, we have $f^{-1}(Y_{sing})\subset U= f^{-1}(V)$. 
Now, we need only show how to deform $f|_{U} \rel \p U$ equivariantly to a $G$-isovariant map, $f'\co U\to Y$. 

But if $q\in Y_{sing}$ and $g\in G$ and $gq\neq q$, then we have $U_q\cap U_{gq}=\varnothing$.
Therefore, if we choose one point $q$ from each $G$-orbit in $ Y_{sing}$ it is sufficient to show how to deform $f|_{U_q} \rel \p U_q$    to a $G_q$-isovariant map $f'_q\co  U_q\to Y$.
Here, $G_q := \{ g \in G \mid gq=q \}$ denotes the isotropy group of $q$.

Fix $q\in Y_{sing}$. 
Let $Y_{(q)}= (Y-Y^G)\cup\{q\}$.   
Recall $f|_{U_q} = e_1\circ \lambda_q\circ f|_{U_q}$, where $\lambda_q\co V_q\to  t_q Y_{(q)}$, and $e_1\co  t_q Y_{(q)} \to Y$ is equivariant.
By \fullref{obstr}, there is a $G_q$-isovariant map $F_q\co U_q  \to t_qY_{(q)} $ for which $(\lambda_q\circ f|_{U_q}) \simeq_{G_q} F_q \rel \p U_q$. 
Define $f'_q: =e_1\circ F_q$. Then $f'_q$ is $G_q$-isovariant and $f|_{U_q}\simeq_{G_q} f'_q $, as required.
\end{proof}

\begin{proof}[Proof of \fullref{equi-iso}(2)]
The argument is entirely similar to that for the proof of \fullref{equi-iso}(1).
Realize the homotopy from $f_0$ to $f_1$ by a $G$-map
\[
(F, f_0\sqcup f_1) ~:~ X\times(I,\p I)  \longrightarrow Y\times(I,\p I).
\] 

\emph{Step 1:} We find 1-neighborhoods $U, V$ of $X_{sing}\times I, Y_{sing}\times I$ so that $F^{-1}(V)=U$.

Choose a tamed neighborhood $(W,h)$ of $Y_{sing}$ in $Y$. Let $W_q$ be the component of  $W$ containing $q$, for each $q\in Y_{sing}$. Since $\dim(Y\times I)\geq 5$, by Siebenmann's thesis~\cite{Siebenmann_thesis} again, we can choose  a 1-neighborhood $V$ of $Y_{sing}\times I$ in $W\times I$, such that $V\cap (Y\times\p I)$ is a $0$-neighborhood of $Y_{sing}\times \p I$ in $Y\times \p I$.

Let $\p_0V$ be the frontier of $V$ in $Y\times I$. Then $\p_0V$ is a codimension $0$-submanifold of $\p V$, and $\p\p_0V=\p_0V\cap (Y\times \p I)$.

Make  $F$ transverse to $\p_0 V$ after a homotopy which is isovariant on $X\times \p I$. Let $U=F^{-1}(V)  $ and let    $U_0=f^{-1}(\p_0 V), $   the  frontier of $U$ in $X\times I$. Then $\p_0U$ is a manifold with boundary and $\p \p_0U \subset  X\times \p I$. Also $(\p_0U,\p \p_0U)$ is bicollared in $(X\times I,\p (X\times I))$.

Proceed as in the proof of \fullref{equi-iso}(1) to make $U$ a $0$-neighborhood of $X_{sing}\times I$. As before let $U_q$ be the component of $U$ containing $(F^G)^{-1}(q\times I)$, and $\p_0U_q= \p_0U\cap U_q$. 

We plan to make $\p_0U_q$ simply connected for each $q$. We repeat the ``innermost circles" argument of 
Browder~\cite{Browder_MxR} doing handle exchanges along $\p_0U_q$ using $2$-handles to kill off the finitely many generators of each $\pi_1(\p_0U_q)$. We do one $q$ at a time, choosing one $q$ from each $G$ orbit of $Y_{sing}$.  In the end we get  a new $F\co X\times I \to Y\times I$ with its new $U$  for which $\p_0U_q$ is $1$-connected for each $q\in Y_{sing}$.  This implies that $\pi_1(X\times I)=\pi_1(U_q)*\pi_1(X\times I -\Int(U_q))$, and $\incl_*\co  \pi_1(U_q)\to \pi_1(X\times I)$ is injective. But $V_q$ is simply connected and therefore $U_q$ is simply connected by the diagram
\[\xymatrix{
\pi_1(U_q) \ar[r]^-{(f|_{U_q })_*} \ar[d]_-{\incl_*} & \pi_1(V_q)=\{1\} \ar[d]^-{\incl_*} \\
\pi_1(X\times I) \ar[r]^-{F_*}_-{\cong} & \pi_1(Y\times I).
}\]

Therefore $U$ is a $1$-neighborhood of  $X_{sing}\times I$ in $X\times I$, and  $U=F^{-1}(V)$.

\emph{Step 2:} In this step we show $f_0$ is isovariantly homotopic to $f_1$.

Let $H=p_1\circ F\co  X\times I\to Y$. Note $H$ is a $G$-homotopy from $f_0$ to $f_1$, with
\[
H^{-1}(Y^G)\subset U\subset H^{-1}(W);  \quad H(U_q)\subset W_q\quad \forall \; q\in Y^G.
\]

As in the proof of  \fullref{equi-iso}(1), for each $q\in Y^G$, $\lambda_q\circ H|_{U_q}\co  U_q\to t_qY_{(q)}$ is homotopic $\rel \p U_q$ to an isovariant map, by \fullref{obstr}.
Therefore $F$ is homotopic $\rel X\times I - \Int(U)$, to an isovariant map $H'\co  X\times I\to Y$, which serves as the required isovariant homotopy from $f_0$ to $f_1$. 
\end{proof}

Observe that we did note use any end theorems in the proof above.  

%------------------------------------------------------------------------------
%------------------------------------------------------------------------------
\section{Quinn--Ranicki = Davis--L\"uck in the case of free actions}\label{sec:QR=DL}

In this section we identify the Quinn--Ranicki assembly map with a map in equivariant homology  in the case of a group acting freely on a CW-complex, where all components are simply connected.  (The connected components may be permuted by the group action.)   For a connected CW-complex this follows from the characterization of assembly maps due to Weiss--Williams \cite{WW}, but in our case we must use the equivariant characterization of assembly maps given in \cite[Section~6]{DL1}.   Hambleton and Pedersen  generalized the work of Weiss and Williams in a different direction. Unfortunately, Davis--L\"uck \cite{DL1} did not connect the map in equivariant homology with the Farrell--Jones Conjecture.   This was remedied in \cite[Corollary~10.2]{HambletonPedersen}, which  identified the stratified assembly map used in the original formulation of the Farrell--Jones Conjecture \cite{FJ_isoconjs} with the map in equivariant homology induced by $E_\vc G \to E_\all G$.  Thus their work~\cite{HambletonPedersen} applied 
to the contractible case, while ours applies to the free case.

Let $\Ho \Spectra$ be the \emph{homotopy category}, given by formally inverting weak homotopy equivalences. There is a \emph{localization functor} $\Ho\co \Spectra \to \Ho \Spectra$ sending weak equivalences to isomorphisms, and this functor is initial with respect to all such functors from $\Spectra$.  The functor $\Ho$ is a bijection on objects. Homotopy groups $\pi_i \co  \Spectra \to \ab$ factor through the functor $\Ho$.  Let $\mcC$ be a category.  A \emph{$\mcC$-spectrum} is a functor from $\mcC$ to $\Spectra$, a \emph{map of $\mcC$-spectra} is a natural transformation, and a \emph{weak equivalence of $\mcC$-spectra} is a map of $\mcC$-spectra $\bE \to \bF$ which induces a weak homotopy equivalence of spectra $\bE(c) \to \bF(c)$ for all objects $c$ in $\mcC$.  
There is a localization functor $\Ho \co  {\mcC}\text{-}\Spectra \to \Ho {\mcC}\text{-}\Spectra$. A key property is that  if $\bE$ and $\bF$ are $\mcC$-spectra which become isomorphic in $\Ho {\mcC}\text{-}\Spectra$, then there is a $\mcC$-spectrum $\bG$ and weak equivalences $\bE \leftarrow \bG \to  \bF$.  

For a groupoid $\mcG$, let $\bL(\mcG)$ be the corresponding $L$-spectrum, as in \cite[Section~2]{DL1}.  This is a functor from the category of groupoids to the category of spectra which satisfies the additional property that an equivalence  $F  \co  \mcG \to \mcG'$ of groupoids induces a weak  equivalence $\bL(F) \co   \bL(\mcG) \to \bL(\mcG')$ of spectra.  
  
Ranicki, motivated by earlier geometric work of Quinn, defined the \emph{assembly map} \cite[Chapter~14]{Ranicki_TOP}, a natural transformation of functors from $\topo$ to $\Spectra$:
\[
\bA(Z) \co  \bH(Z; \bL(1)) \longra \bL(\Pi_1 Z)
\]
where $\Pi_1 Z$ is the fundamental groupoid of $Z$. When $Z$ is a point, the assembly map is a weak equivalence. 
The \emph{algebraic structure spectrum $\mcS^{per}(Z)$} is defined to be the homotopy cofiber of $\bA(Z)$. Its homotopy groups $\mcS^{per}_*(Z) := \pi_*\mcS^{per}(Z)$ are the algebraic structure groups used in \fullref{sec:calculation}; one can do this for pairs also.
 
Fix a group $G$. Consider the orbit category $\Or(G)$ and the $\Or(G)$-spectrum
\[
\ul{\bL} \co  \Or(G) \longra \Spectra ~;~  G/H \longmapsto \bL(\overline{G/H})
\]
where $\overline{G/H}$ is the groupoid associated to the $G$-set $G/H$.  
For a $G$-CW-complex $X$, consider the spectrum 
\[
\bH^G(X; \ul{\bL}) ~:=~ \map_G(-,X)_+ \wedge_{\Or(G)} \ul{\bL}(-).
\]
Then, by definition, $H_*^G(X;\ul{\bL}) = \pi_*\bH^G(X; \ul{\bL})$.
 
Write $G\cw$ for the category whose objects are $G$-CW-complexes and whose morphisms are cellular $G$-maps.  (Actually, for set-theoretic reasons we need to restrict our $G$-CW-complexes to a fixed universe; for our purposes it will suffice to assume that the underlying space of each CW-complex is embedded in $\R^{\infty}$.)  A \emph{$(G,\mcF)$-CW-complex} is a $G$-CW-complex with isotropy in a family $\mcF$.  
Let  $(G,\mcF)\cw$ be the full  subcategory of $G\cw$ whose objects are $(G,\mcF)$-CW-complexes.   Let $\Or(G,\mcF)$ be the full subcategory of $(G,\mcF)\cw$ whose objects are the discrete $G$-spaces $G/H$ with $H \in \mcF$.
The symbol $1$ will denote both the trivial group and the family of subgroups of $G$ consisting of the trivial group.  Let $\mathrm{sc}(G,1)\cw$ be the full subcategory of $(G,1)\cw$ whose objects are free $G$-CW-complexes all of whose components are simply connected.

Let $X$ be a free $G$-CW complex.   Let $\Pi_0X$ be the $G$-set of path components of $X$.
Here is the main theorem of this appendix. 

\begin{thm}\label{QR=DL}
There is a commutative diagram in $\Ho\,(G,1)\cw\text{-}\Spectra$:
\[\xymatrix{
\bH(X/G; \bL(1)) \ar[r] \ar[d] & \bL(\Pi_1(X/G)) \ar[d]\\
\bH^G(X; \ul{\bL}) \ar[r] & \bH^G(\Pi_0 X; \ul{\bL}).
}\]
\begin{enumerate}
\item  The top map is the assembly map $\bA(X/G)$ and is a map of $(G,1)\cw\text{-}\Spectra$.
\item  The bottom map is induced by the $G$-map $X \to \Pi_0X$ and is a map of $(G,1)\cw\text{-}\Spectra$.
\item  The right map is the composite of the formal inverse of the weak equivalence of $(G,1)\cw$-spectra  $\bL(\Pi_1(EG \times_G X)) \to \bL(\Pi_1(X/G))$ and the map of $(G,1)\cw$-spectra $\bL(\Pi_1(EG \times_G X)) \to \bH^G(\Pi_0X; \ul{\bL})$ defined in \fullref{lem:groupoids}.  This map is a weak equivalence when restricted to $ \mathrm{sc}(G,1)\cw\text{-}\Spectra$.
\item  The left map is an isomorphism in $\Ho\,(G,1)\cw\text{-}\Spectra$.
\end{enumerate}
\end{thm}

The proof of the theorem is quite formal and applies more generally.  What is needed is a functor from groupoids to spectra which takes equivalences of groupoids to weak equivalences of spectra and an assembly map which is a weak equivalence when $X$ is a point.  So, for example, our theorem applies equally well to $K$-theory.  See \fullref{non-orientable} below for the modifications necessary for the $L$-theory non-orientable case.

Let $G$ be a discrete group. Let $S$ be a $G$-set. Define the \emph{action groupoid} $\overline{S}$ as the category whose object set is $S$, and whose morphisms from $s$ to $t$ are triples $(t,g,s)$ such that $t=gs$, and whose composition law is $(t,g,s) \circ (s,f,r) = (t,gf,r)$.  Define a functor 
\[
\bL^G\co  (G,1)\cw\ \longrightarrow \Spectra ~;~ X \longmapsto \bL(\Pi_1(EG \times_G X)).
\]

The next lemma relates $\bL^G$ to the above functor $\ul{\bL}\co  \Or(G) \to \Spectra$.

\begin{lem}\label{lem:groupoids}
Let $G$ be a discrete group.

\begin{enumerate}
\item\label{spectra_homeo}
For a discrete $G$-set $S$, there is a homeomorphism of $spectra$
\[
\bH^G(S; \ul{\bL}) \cong \bL(\overline{S}),
\]
natural in $S$.

\item\label{Phi}
For a free $G$-CW-complex $X$, there is a map of groupoids 
\[
\Phi(X) \co  \Pi_1(EG \times_G X) \longra \overline{\Pi_0X} 
\]
which is an equivalence of groupoids when all the components of $X$ are simply connected.  Furthermore,  $\Phi$ is natural in $X$; that is, $\Phi(-)$ is a map of $(G,1)\cw$-groupoids.  

\item\label{Lambda}
There is a map of $(G,1)\cw$-spectra
\[
\Lambda(X) \co  \bL^G(X) \longra \bH^G(\Pi_0 X; \ul{\bL})
\]
whose restriction to $\mathrm{sc}(G,1)\cw$ is a weak equivalence.
\end{enumerate}
\end{lem}

\begin{proof}
\eqref{spectra_homeo}  The homeomorphism is given by
\[
\bH^G(S ; \ul{\bL}) \longra \bL(\overline{S}) ~;~ [(f,x) \in \map_G(G/K,S)_+ \wedge L(\overline{G/K})_n] \longmapsto L(\overline{f})_n(x) \in L(\overline{S})_n.
\]
If $S$ is an orbit $G/K$, then the inverse is given by $x \in L(\overline{G/K})_n \mapsto [(\id,x) \in \map_G(G/K,G/K)_+ \wedge L(\overline{G/K})_n]$.  The case of a general $G$-set 
follows since both $\bH^G(-; \ul{\bL})$ and $\bL(-)$ convert disjoint unions to one-point unions of spectra.  

\eqref{Phi}  We first need some notation.  For a subset $A$ of a topological space $Y$, let $\Pi_1(Y,A)$ be the full subcategory of the fundamental groupoid $\Pi_1Y$ whose objects are points in $A$.  If $\Pi_0A \to \Pi_0Y$ is onto, then there is an equivalence of groupoids $\Pi_1 Y \to \Pi_1(Y,A)$ whose definition depends on a choice of a  path from $y$ to a point in $A$ for every $y \in Y$.

Let $p \co  EG \times X \to EG \times_G X$ be the quotient map.  We will define $\Phi(X)$ as a composite
\[\xymatrix{
\Pi_1(EG \times_G X) \ar[r]^-{\Theta(X)} & \Pi_1(EG \times_G X, p(\{e_0\}\times X)) \ar[r]^-{\Psi(X)} & \overline{\Pi_0X}.
}\]

We first define $\Theta(X)$ by making choices in the universal space $EG$.  Choose a point $e_0 \in EG$.  For each $e \in EG$, choose a path $\sigma_e \co  I \to EG$ from $e$ to $e_0$, choosing the paths so that, for all $g \in G$ and $t \in I$, $g(\sigma_e(t)) = \sigma_{ge}(t)$.  This can be accomplished by choosing a set-theoretic section $s \co  BG \to EG$ of the covering map, and defining the remaining $\sigma_e$ by equivariance.  Then for $p(e,x) \in EG \times_G X$, define the path $\theta_{p(e,x)}(t) := p(\sigma_e(t),x)$.  This path is independent of the choice of representative of $p(e,x)$.  These paths give the equivalence of groupoids $\Theta(X)$, natural in $X$.  

We now define $\Psi(X)$ using the fact that $p$ is a covering map.   On objects, define $\Psi(X)(p(e_0,x)) := C(x) \in \Pi_0X$, where $C(x)$ is the path component of $x$ in $X$.  For  a morphism represented by a path $\alpha \co  I \to EG \times_G X$ with $\alpha(0) = p(e_0,x)$ and $\alpha(1) = p(e_0,y)$, let $\widetilde{\alpha} \co  I \to EG \times X$ be the lift of $\alpha$ starting at $(e_0,x)$.  Then $\widetilde{\alpha}(1) = (ge_0, gy)$ some $g \in G$.  Then define $\Psi(X)[\alpha] := (C(x), g, C(y))$.  We leave the geometric details of verifying that this is a functor to the reader, but note that we follow that convention that a path $\alpha$ determines a morphism from $\alpha(1)$ to $\alpha(0)$ in the fundamental groupoid.

Suppose all the components of $X$ are simply connected. We now show that $\Psi(X)$ is an equivalence of groupoids. Choose a base point for each component of $X$.  Define a functor $\overline{\Pi_0X} \to \Pi_1(EG \times_G X, p(\{e_0\} \times X))$ on objects by sending $C(x)$ to $p(e_0,x)$, and on morphisms by sending $(C(x), g, C(y))$ to $[p \circ \widetilde{\alpha}]$ where $\widetilde{\alpha} \co  I \to EG \times X$ is a path from $(e_0,x)$ to $(ge_0, gy)$. This $\widetilde{\alpha}$ is unique up to homotopy rel endpoints since $X$ is simply connected. This ends the proof of \eqref{Phi}.

\eqref{Lambda}  Define $\Lambda(X)$ as the composite of $\bL(\Phi(X))$ and the isomorphism from $\eqref{spectra_homeo}$. 
\end{proof}

\begin{rem}\label{rem:excisive}
We next recast the axiomatic approach of \cite[Section~6]{DL1}.  Our terminology is self-consistent but does not precisely match that of \cite{DL1}; in particular we drop the adverb ``weakly.''  A functor $\bE \co  (G,\mcF)\cw \to \Spectra$ is  \emph{homotopy invariant} if any homotopy equivalence induces a weak equivalence of spectra.  A functor $\bE \co  (G,\mcF)\cw \to \Spectra$ is \emph{excisive} if $\bE(-)$ and $\bH^G(-; \bE|_{\Or(G,\mcF)})$ are isomorphic objects in $\Ho\,(G,\mcF)\cw\text{-}\Spectra$.  This is equivalent to the notion of weakly $\mcF$-excisive given in \cite{DL1}.\footnote{Indeed \cite[Theorem~6.3(1,3)]{DL1} implies that if $\bE$ is weakly $\mcF$-excisive in the sense of \cite{DL1}, then there is a $(G,\mcF)$-spectrum  $\bE^{\%}$ and weak equivalences  $\bE \leftarrow \bE^{\%} \rightarrow \bH^G(-; \bE|_{\Or(G,\mcF)})$.
Conversely, if $\bE(-)$ and $\bH^G(-; \bE|_{\Or(G,\mcF)})$ are isomorphic objects in $\Ho\,(G,\mcF)\cw\text{-}\Spectra$, then there are weak equivalences $\bE \leftarrow \bF \rightarrow \bH^G(-; \bE|_{\Or(G,\mcF)})$ for some $(G,\mcF)\cw$-spectrum $\bF$.  But \cite[Theorem~6.3(1)]{DL1} shows that $\bH^G(-; \bE|_{\Or(G,\mcF)})$ is weakly $\mcF$-excisive, and hence so is any weakly equivalent $(G,\mcF)\cw$-spectrum.}  By \cite[Theorem~6.3(2)]{DL1}, a map $\bT \co  \bE \to \bF$ of excisive $(G,\mcF)$-CW-spectra is a weak equivalence if and only if $\bT(G/H) \co  \bE(G/H) \to \bF(G/H)$ is a weak equivalence of spectra for all $H \in \mcF$.
An \emph{excisive approximation} of a homotopy invariant functor $\bE\co  (G,\mcF)\cw \to \Spectra$ is a map $\bT\co  \bE' \to \bE$ of $(G,\mcF)$-CW-spectra such that $\bE'$ is excisive and $\bT(G/H)$ is a weak equivalence for all orbits $G/H$ with $H \in \mcF$.

We next assert existence and uniqueness of excisive approximations.  Theorem~6.3(2) of \cite{DL1} constructs a specific excisive approximation $\bE^{\%} \to \bE$ which is functorial in  $\bE$.  Excisive approximations are unique in the sense that, given any two excisive approximations $\bT' \co  \bE' \to \bE$ and $\bT'' \co  \bE'' \to \bE$, there is an isomorphism $\bS \co \bE' \to \bE''$ in $\Ho\,(G,\mcF)\cw\text{-}\Spectra$ such that $\bT' = \bT'' \circ \bS$.  Indeed, to verify that $\bS$ exists, it suffices to compare any excisive approximation $\bT' \co  \bE' \to \bE$ is equivalent to the functorial  excisive approximation $\bT^{\%} \co  \bE^{\%}  \to \bE$.  Consider the commutative diagram in $(G,\mcF)\cw\text{-}\Spectra$:
\[\xymatrix{
\bE'^{\%} \ar[r] \ar[d] & \bE' \ar[d]\\
\bE^{\%} \ar[r] & \bE.
}\]
Since the left and top maps are both weak equivalences,  we obtain an isomorphism in the homotopy category $\bS := (\bE' \leftarrow \bE'^\% \rightarrow \bE^\%)$ with $\bT' = \bT^{\%} \circ \bS$, as desired.
\end{rem}

\begin{proof}[Proof of \fullref{QR=DL}]
The theorem will be proven by concatenating three commutative squares.   The first is a commutative diagram in $\Ho\,(1,1)\cw\text{-}\Spectra$, which we will apply below in the case $Z = X/G$.
\[\xymatrix{
\bH(Z; \bL(1)) \ar[r] \ar[d] & \bL(\Pi_1 Z) \ar[d]\\
(\bL^1)^{\%}(Z) \ar[r] & \bL^1(Z)
}\]
The right map is the identity.  The top map is the assembly map $\bA(Z)$.  The left map exists in the homotopy category (see end of \fullref{rem:excisive}) and is an isomorphism since both horizontal maps are $(1,1)$-excisive approximations of $\bL(\Pi_1 Z) = \bL^1(Z)$.

Next comes a commutative diagram in $\Ho\,(G,1)\cw\text{-}\Spectra$.  
\[\xymatrix{
(\bL^1)^{\%}(X/G) \ar[r] & \bL^1(X/G)\\
(\bL^G)^{\%}(X) \ar[r] \ar[u] & \bL^G(X) \ar[u]
}\]
The right map is induced by the homotopy equivalence $EG \times_G X \to X/G$, inducing an equivalence of fundamental groupoids, hence a weak equivalence of spectra.  The left map exists and is an isomorphism, since the top map and the composite of the bottom and right maps are $(G,1)$-excisive approximations of $\bL^G(X)$.

Our final commutative diagram is in $\Ho\,(G,1)\cw\text{-}\Spectra$.
\[\xymatrix{
(\bL^G)^{\%}(X) \ar[r] \ar[d] & \bL^G(X) \ar[d]\\
\bH^G(X; \ul{\bL}) \ar[r] & \bH^G(\Pi_0X; \ul{\bL})
}\]
The top map  is the functorial excisive approximation of $\bL^G(X)$. The bottom map is induced by the $G$-map $X \to \Pi_0X$ and is an excisive approximation of $ \bH^G(\Pi_0 (-); \ul{\bL})$.   The right map is defined in \fullref{lem:groupoids}\eqref{Lambda} and is an isomorphism when restricted to $\Ho \mathrm{sc}(G,1)\cw\text{-}\Spectra$.  Functoriality gives a map $(\bL^G)^{\%}(X)  \to \bH^G(\Pi_0 (-); \ul{\bL})^{\%}(X)$, and the bottom map is an excisive approximation implies $\bH^G(\Pi_0 (-); \ul{\bL})^{\%}(X) \to  \bH^G(X; \ul{\bL})$; define the left map as the composite.  Since the left map is a map of excisive functors and  is a homotopy equivalence when $X = G/1$, the left map is an isomorphism in $\Ho\,(G,1)\cw\text{-}\Spectra$.
\end{proof}

\begin{rem}\label{non-orientable}
We next indicate the modifications needed for the statement and proof of \fullref{QR=DL} in the non-orientable case.  A \emph{groupoid with orientation character} $\mcG^w$ is a groupoid $\mcG$ together with a  functor $w \co  \mcG \to \{ \pm 1 \}$, where $\{\pm 1\}$ is the category with a single object and two morphisms $\{+1,-1\}$ where $-1 \circ -1 = +1$.  A map of groupoids
with orientation character is a map of groupoids which preserves the orientation character.  Let $\matheurm{GWOC}$ denote the category of groupoids with orientation character. There is an $L$-theory functor $\bL \co  \matheurm{GWOC} \to \Spectra$. (The definition in \cite{DL1} can be easily modified to cover this case, see also \cite{BL_crossed}.)  Two maps $F_0,F_1 \co  \mcG^{w} \to  \mcG'^{w'}$ of groupoids with orientation character are \emph{equivalent} if there is a natural transformation which is orientation preserving in the sense that $w'(F_0(x) \to F_1(x)) = +1$ for all objects $x$ of $\mcG$.  A map $F  \co  \mcG^{w} \to  \mcG'^{w'}$ is an \emph{equivalence of groupoids with orientation characters} if there is a map $F' \co  \mcG'^{w'} \to  \mcG^{w}$ so that both composites $F \circ F'$ and $F' \circ F$ are equivalent to the respective identity.  An equivalence of groupoids with orientation characters gives an weak equivalence of $L$-spectra.

Now suppose $G$ is a group with orientation character $w \co  G \to \{\pm 1\}$.  Following forthcoming work of Davis and Lindenstrauss, we discuss two related groupoids with orientation character.  First, if $S$ is a $G$-set,  give the action groupoid $\overline{S}$ the orientation character $(t,g,s) \mapsto w(g)$.  This gives a functor
 $\Or(G) \to  \matheurm{GWOC}$ defined on objects by $G/H \mapsto (\overline{G/H})^w$ and hence a functor:
 \[
\ul{\bL} \co  \Or(G) \longrightarrow \Spectra ~;~ G/H \longmapsto \bL((\overline{G/H})^w).
\]

Suppose $\phi \co  \widehat Y \to Y$ is a double cover.    Define the \emph{fundamental groupoid with orientation character} $\Pi_1^w(Y)$, as follows.  The objects are the points of $\widehat Y$.  A morphism from $\hat{y}$ to $\hat{y}'$ is a path $\alpha$ from $w(\hat{y}')$ to $w(\hat{y})$.  A morphism is assigned $+1$ if the unique lift of $\alpha$ starting at $\hat{y}'$ ends at $\hat{y}$; otherwise assign the morphism $-1$. 

Recall $G$ is a group with orientation character $w$.  Given  a free $G$-CW-complex $X$, let 
$w \co  EG \times_{\Ker(w)} X \to EG \times_G X$ be the corresponding double cover.
Thus, for a fixed $(G,w)$, there is a functor $\bL^G$ defined by
\[
\bL^G\co  (G,1)\cw\ \longrightarrow \Spectra ~;~ X \longmapsto \bL(\Pi_1^w(EG \times_G X)).
\]

Then, after modifying $\bL$, $\bL^G$ and $\overline{S}$ as indicated above, the statement and proof of \fullref{lem:groupoids} remain valid.  The same is true for \fullref{QR=DL} after accounting for Ranicki's version of the assembly map in the non-orientable case \cite[App.~A]{Ranicki_TOP}.
\end{rem}

%------------------------------------------------------------------------------
%------------------------------------------------------------------------------
\bibliographystyle{gtart}
\bibliography{cdk}

\end{document}